\documentstyle[11pt,a4wide]{article}
\title{\bf A Complete Classification Of The Admissible  Representations
 Of Infinite-Dimensional Classical Matrix Groups}
\author{N. I. Nessonov}
\date{}
\begin{document}
\maketitle

  In the present paper a complete description is obtained for the class of
 admissible (by the terminology of  G. I. Ol'shansky) unitary
 representations of infinite-dimensional analogous  classical matrix groups
 $GL(\infty), Sp(2\infty), O(2\infty).$

Henceforth these objects we will imagine as matrix of operators that acts in
the complex Hilbert space $H$.  In a case $GL(\infty)$ suppose, that basis
$\lbrace e_1,e_2,\ldots,e_n,\ldots \rbrace$ of $H$ is numerated by the
elements from  ${\bf N}$, and $e_n$ is an infinite row with the unique
nonzero  $\; {\it n}$-th coordinate that equals to unit.

Let $B(H)$ be a set of bounded operators in $H$, $GL(H) = \lbrace g\in B(H)
: $ there  exist $g^{-1} \in B(H) \rbrace, GL(n) = \lbrace g \in GL(H) :
ge_i = g^*e_i = e_i $ for all $i > n \rbrace$, and $GL(\infty)$ is an
inductive limit of the groups $GL(n)\;( n = 1, 2, 3,... )$.

To definite $Sp(2\infty), O(2\infty) $ let's consider in $H$ orthonormal
 basis \newline $\lbrace \ldots, e_{-n}, \ldots, e_{-2}, e_{-1}, e_1, e_2, \ldots,
e_n,\ldots \rbrace$  and represent $H$ in the form $H = H_-\oplus H_+$,
where $H_-$ and
$H_+$ are generated by basis elements with the negative and positive numbers
respectively. We put $GL(2n) = \lbrace g \in GL(H) : ge_i = g^*e_i = e_i,$
if $\vert i \vert > n \rbrace$ and define operators $s_+$ and $s_-$
in $H$ by the formulae :$$s_+e_i = e_{-i}, s_-e_i = (sign\, \ \
i)e_{-i}.$$

Let $GL(2 \infty)$ be an inductive limit of $GL(2n)$, $Sp(2\infty)$
$(O(2\infty)) = \lbrace g\in GL(2\infty) : s_-g^ts_-^{-1} = g^{-1}
(s_+g^ts_+ = g^{-1}) \rbrace$, where $g^t$ is a transposed matrix with
respect to $g$.

We will denote  by $G$ one of the groups : $GL(\infty)$, $Sp(2\infty)$ or
$O(2\infty)$. If a more specific setting will be necessary, we'll specially
notice it.  We denote by $U(G)$ the unitary subgroup of group $G$ and remind
the definition of the admissible representation.  \medskip

{\bf{Definition}} ( see \cite{Olsh1} and \cite{Olsh2} ).
  {\it{  Let $G(n, \infty) = \lbrace g \in G :
 ge_i = g^*e_i = e_i$ for all $ i : \vert i \vert < n \rbrace$.
  Factor-representation ${\Pi} $ of group  $G$, that acts in
  the Hilbert space
$H_{_\Pi},$ is called
 {\sl an admissible} one, if there exists $n \in { \bf N} $ and nonzero
  vector $ \xi \in H_{_{{\Pi}}}$
 with the property: ${\Pi}(u) \xi = \xi$ for all
  $u \in U(G(n, \infty))$.}}
\medskip

In some chapters of this paper it's better to use the equivalent

{\bf{Definition}} of an admissible representation (see \cite{Olsh1}).
  {\it{ The representation ${\Pi}$
 is called  an {\sl{admissible}} one,
  if the set $${\bigcup_{n=1}^{\infty}\lbrace
\eta \in H_{_{\Pi}} :  {\Pi}(u) \eta = \eta\;
\mbox{for all}\;  u \in U(G(n, \infty)\rbrace}$$
  is dense in $H_{_{\Pi}}$}}.
\medskip

Let's denote by $G_n^I$  the subgroup $GL(\infty)$, which consists
  of matrices $\pmatrix{I_n & 0\cr \ast & \ast}$,
where $I_n\; \mbox{is a unit}\;n \times n-$  matrix.
  When $n=0,\; \mbox{then}\;\; G_0^I = GL(\infty)$.

Let's give a construction of the representations of group $G_n^I$,
  that will be very important  in our  following reasoning.

Let $\Lambda_m$ be a set of all complex matrices of $m$ rows and infinite
number of columns, $\nu_m-$  the Gaussian  measure on $\Lambda_m$
 with the unit covariant
operator, $A$ an $m \times m$ selfadjoint  matrix, $z$ a matrix of $n$
rows and $m$ columns ($ z - n \times m-$ matrix).
\medskip

 Let's define the representation ${\Pi}_{_{\scriptstyle{Az}}}$
  of group $G_n^I$ ( ${\Pi}_{_{\scriptstyle A}}$ of group $GL(\infty)=
G_0^I)$ in $L^2(\Lambda_m,\nu_m)$ by

$$({\Pi}_{_{\scriptstyle{Az}}}{\bigg (}\pmatrix{I_n & 0 \cr 0 & g}
 {\bigg )}\eta)(\lambda)= \mid {\rm det}  g \mid ^{i \beta}
{{\hat{\alpha}_{_A}} } (\lambda,g)\eta(\lambda g), $$
where$${\hat{\alpha}_{_A}}(\lambda,g) = \exp{\bigg\lbrace} \frac {-1}{2}
Tr \lbrack \lambda(gg^{\ast}-1)\lambda^{\ast}-
2iA\lambda(gg^{\ast}-1)\lambda^{\ast} \rbrack {\bigg\rbrace}, $$
$I_n\;$ is $\; n\times n$ --unit matrix,  $\beta\,-$ real number;
$$({\Pi}_{_{\scriptstyle{Az}}}{\bigg (}\pmatrix{I_n & 0 \cr h & I}
 {\bigg )}\eta)(\lambda)=
\exp \lbrack i \, \  \Re \, \  Tr (z \lambda h ) \rbrack
\eta (\lambda) $$
$$( \eta \in L^2 ( \Lambda_m, \nu_m)).$$

If $M\,$ is a set of operators in the Hilbert space, and $M^{ \prime}$
is a commutant of $M$, then we can write down
\medskip

{\bf{Theorem 0.1.}}{\it{ Von Neumann algebra ${({\Pi}_{_{\scriptstyle{Az}}}
  (G_n^I))}^{ \prime} $ is generated by
the operators $\tau(u),$ where  $u \in \lbrace v \in U(m): \, \ \   vA=Av
\, \ \   and \, \ \   zv=z \rbrace=U(m,A,z), $
that acts in $L^2(\Lambda_m,\nu_m)$ by:
 $ (\tau(u)\eta)(\lambda)=\eta(u^{\ast}\lambda).$}}

\medskip
At first this fact was announced in \cite{Olsh3}  for group $GL(\infty).$

In \cite{N1} and \cite{N2}  the complete description of spherical
  representations of
$G_n^I\;$ was obtained. It's maintained  in the following  statement:
\medskip

{\bf{Theorem 0.2.}}{\it Let ${\Pi}$ be a factor-representation
  of $G_n^I,$
that acts in a Hilbert space $H_{_{ \Pi}}.$ If there exists in $H_{_{ \Pi}}$
a nonzero vector $\xi$, that fixed by the operators ${\Pi}(U(G_n^I)),$
where $U(G_n^I)$ is the unitary subgroup $G_n^I$, then ${\Pi}$ is multiple
 (for some
m,A,z) by restriction of ${\Pi}_{_{\scriptstyle{Az}}}$
  to subspace $\lbrace \eta \in L^2(\Lambda_m,\nu_m) :
\tau(u)\eta=\eta \, \ \  $ for all $\, \ \   u \in U(m,A,z) \rbrace.$}
\medskip

If $\rho$ is an irreducible representation of $U(m,A,z)$,
${{\rho}_{_{kl}}} ({1} \le k,l \le \dim {\rho})$ is its matrix element,
  then operator
\newline $P_{k\rho}= \dim\rho\int \limits_{U(m,A,z)}
  {\bar\rho}_{_{kk}} \tau(u) du$
is a minimal orthoprojection from $({\Pi}_{_{\scriptstyle{Az}}}
 (G_n^I))^ \prime. $

The  main result of the paper in a case of group\, \ \    $G_n^I$   \ \ \
 is a following
\medskip

{\bf{Theorem 0.3.}}{\it Let ${\Pi}$ be an admissible factor-
 representation of group
$G_n^I.$ Then there exist    $\; m,A,z,\rho$\, \ \   such
 that ${\Pi}$ is multiple by
restriction of ${\Pi}_{_{\scriptstyle{Az}}}$ to
 $P_{_{\scriptstyle{k\rho}}}L^2(\Lambda_m,\nu_m).$  }

\medskip
Classification result (see Theorem 6.15) for the groups
 $Sp(2\infty)$
and $\; O(2\infty)$ is obtained by the same way.

Namely, in \S 1  there was built the set of unitary representations
  (reducible ones)
${\Pi}_{_{\scriptstyle{A}}}  $ of these groups, that have almost
 the same meaning,
as representation    $\; {\Pi}_{_{\scriptstyle{Az}}}\, $
 for $ G_n^I$.
In Proposition 1.4  the form of their decomposition into  irreducible
 components
is obtained.
\medskip

The main result of the paper in a case of groups $Sp(2\infty)$
and  $\; O(2\infty)$ contains the following theorem, which is
 proved in \S 6 (see in greater detail theorem 6.15).
\medskip

{\bf{Theorem 0.4.}}{\it { Let $\Pi$ be an admissible factor--representation
 of group
 $Sp(2\infty)$  or\, \ \   $ O(2\infty)$. Then there exists
$m\times m$-- selfadjoint matrix $A$
such that $\Pi$ is multiple by one of the irreducible components of
 the representation $\Pi_{_{\scriptstyle A}}, $
that was built in propositions 1.2 --- 1.3\, \ \
 {\rm  (see proposition 1.4).}}}
\vskip 16pt
\medskip

Give a brief account  of the logic of our reasoning.
The offered method essentially bases on the following statement,
 that belongs to G.I. Ol'shansky.

\vskip 16pt

{\bf{Theorem 0.5.}}   {\it If $\Pi\;$ is an admissible representation
 of group $G,$
that acts in a Hilbert space $H_{_{ \Pi}},$ $\Pi_{_U}\;$
is its restriction on $U(G),$ then $\Pi_{_U}$
 extends by continuity
to the representation  of semigroup of partial isometrics
 $\tilde U(G),$
which are the limiting points of elements from $U(G)$ concerning to the weak
topology in $B(H)$.} \medskip

All the main classification results (theorems 5.7, 5.10, 6.14, 6.15) follow
 from the structure of the spherical representations of group $GL(\infty)$
 (see \cite{N1}) and group of motions
   $\; G_n^I$  (see\cite{N2}), which we'll identify with
 the subgroup of $\; GL(\infty)\; $,   that consists of matrices
$ \pmatrix{I_n & 0 \cr {\ast} & {\ast}} \, \  ,\;$ where $I_n\;
\mbox{is a unit}\; n\times n-$  matrix, $\ast$ is an arbitrary matrix of
corresponding size (see theorem 2.1).  \medskip

For first let's present the idea of the classification of the admissible
representations of $GL(\infty)$ and $G_n^I$.  \medskip

Let $G$  be the one of these groups. If $G$ =  $GL(\infty),$  then, as
before, $G(p, \infty) = \lbrace g \in G : ge_i = g^*e_i = e_i$ for all $ i :
 i  \leq n \rbrace.$ When $G$ =   $G_n^I$, then we put $G(p, \infty) =
\lbrace g \in G : ge_i =   e_i$ for all $ i :  i  \leq n+p   $ and $ge_i =
g^* e_i $ = $ e_i     $ for all $i$ : $n < i \leq p+n   \rbrace.$

If $\Pi$ is an admissible factor--representation of group $G\, ,$ then for
sufficiently large $p$ in $H_{_{ \Pi}}$ there exists a unit
 $\Pi(U(G(p, \infty))-$ fixed vector $ \xi (p).$

We assume without loss of generality, that  $H_{_{ \Pi}} \, \  =
\, \  [\Pi(G) \xi(p)] \, $,
 where  $ [\Pi(G) \xi(p) ] $ is a closure of linear cover
 of the set $\, \  \{ \Pi(G) \xi(p) \} \, \
(g\in G).$
\medskip

The groups, we consider, have so-called property of asymptotic
 abelianness (see definition 3.1), that allows to define
 for $\Pi$  an important
invariant --an asymptotic spherical function\newline (a. s. f.)
$\varphi_{_{\Pi}}$ (see proposition 3.2, definition 3.3 and
theorem 4.1). It makes it possible for us to define the rang
 $ {\bf r}(\Pi)$
of representation $\Pi.$
\medskip

Henceforth, using the structure of the spherical representation
 of group $G$ and a fact,
 that restriction of $\Pi$ to $G(p, \infty),$ that acts in
 $ \, \  [\Pi(G(p, \infty)) \xi(p)] \, \  \subset H_{_{ \Pi}},$
is irreducible, we prove, that for $p > {\bf r} (\Pi)$  there
 exists a set of the unit

$\Pi(U(G(2(p+1)))-$  fixed vectors $\xi_{_i}^{^U}$ such that

{\em a) in a case when $G = GL(\infty)$   the subspaces
 $ \, \  [\Pi(G_{2(p+1)}^I  \xi_{_i}^{^U} ] \, \  =
 \, \  H_i$

are orthogonal in pairs and $$\bigoplus_i H_i\, \  =
 \, \   [\Pi(G) \xi(p)] \, \  =\, \  H_{_{\Pi}}$$}

 (see lemma 5.3 );

{\em b) when $G = G_n^I$ the subspaces
 $ \, \  [\Pi(G_{2(p+1)+n}^I )
  \xi_{_i}^{^U} ] \, \  = \, \  H_i$

 are orthogonal in pairs and  $$\bigoplus_i H_i \, \  =
 \, \   [\Pi(G) \xi(p)] \, \  = \, \  H_{_{\Pi}}$$}

(see lemma 5.8 ).
\medskip

 That's why for any vector $\eta \in \, \  H_{_{\Pi}}$ in
 the commutant of
$\Pi(G(2(p+1), \infty ) )$ there exists an orthoprojection $P_{_{\eta}}$
 such that for some
natural $i(\eta)$ \, \  $ [\Pi(G(2(p+1),\infty ) )
 P_{_{\eta}}\eta] \subset H_{i(\eta)}.$
Change if we need the orthoprojection $P_{_{\eta}}$ to the lower one, we
suppose, that representation   $\;\; (\Pi,\; G(2(p+1),\infty ),\;
   [\Pi(G(2(p+1),\infty ) ) \eta] )\;$ ({\em a restriction of
 $ \Pi$ to $ G(2(p+1),\infty ),$ that acts in  \, \
 $ [\Pi(G(2(p+1),\infty ) ) \eta]$})
 is multiple by $(\Pi,\; G(2(p+1),\infty ),
 \, \  P_{_{\eta}} H_{_{\Pi}}).$
\medskip

Therefore  $(\Pi,\; G(2(p+1),\infty ) ,  \, \  P_{_{\eta}}
 H_{_{\Pi}})$ is
a restriction of the direct integral of the irreducible representations
 of groups
  $ G_{2(p+1)+n}^I$  to the  $ G(2(p+1),\infty )$
( $ n=0$ corresponds to the case when  $G \, \ =\, \  GL(\infty)$).
Besides that, taking to the consideration statements 5.4 -- 5.5, 5.9
 $ (\Pi,\; G_{2(p+1)+n}^I,\;  H_{i(\eta)})$
 \, \  we can realize in such form  that the restriction
 of every irreducible component
to \, \  $  G(2(p+1),\infty ) \, \  $  is the
 same representation
${\Pi}_{_{\scriptstyle{Az}}}$ ( see lemma 5.5, proposition 5.6 and lemma 5.9 ).
\medskip

Therefore, the irreducible components of representation
 $$(\Pi, \; G(2(p+1),\infty ) ,  \, \  P_{_{\eta}} H_{_{\Pi}})$$
are unitary equivalent to the irreducible components
 of representation  \newline $( {\Pi}_{_{\scriptstyle{Az}}},\;
  G(2(p+1),\infty ),\; L^2(\Lambda_{_{{\bf r}(\Pi)}}, \nu_{_{{\bf
r}(\Pi)}})). $ \medskip

Now let's account an isometry  $\sigma_q^{(n)}$, that acts in  $H$ by
$$\sigma_q^{(n)}(e_i) = e_i\,{\rm when}\, i\leq n\,\mbox{and}\,
\sigma_q^{(n)}(e_i) = e_{i+q}\, {\rm when}\, i> n\,.$$

If element--matrix $g\, =\, \pmatrix{I_{n}&0\cr h_{n}&g_{n}\cr}$\, \
belongs to    $ G,$ then $ \sigma_q^{(n)} g  {( \sigma_q^{(n)})}^*\, \ =\, \
 \pmatrix{I_{n}&0&0\cr 0&0_{q}&0\cr h_{n}&0&g_{n}\cr}$\, \ \newline ($0_q\,
\ $ is $ q\times q $  zero--matrix).  \medskip

Let $g_{\sigma}\,=\, \pmatrix{0&0&0\cr0&I_{q}&0\cr0&0&0\cr}\, +\,
\sigma_q^{(n)} g {( \sigma_q^{(n)})}^*\,$.  Obviously $g_{\sigma} \in
G(q,\infty)$. \medskip

 By theorem 0.5   $\Pi$  extends by continuity on    \, \
 $ \sigma_q^{(n)}$  \, \ , \, \  ${( \sigma_q^{(n)})}^*$
 and for any nonzero vector $ \eta \in H_{_{\Pi}}$
 $(\Pi, G(q,\infty), [\Pi(G(q,\infty))\Pi(\sigma_q^{(n)}) \eta])$
 is multiple by
one of the irreducible components of the representation \hskip 4 pt
 $$( {\Pi}_{_{\scriptstyle{Az}}},
  G(2(p+1),\infty ), L^2(\Lambda_{_{{\bf r}(\Pi)}}, \nu_{_{{\bf
r}(\Pi)}})),$$ when $q=2(p+1).$ Therefore in\, \  $ L^2(\Lambda_{_{{\bf
r}(\Pi)}}, \nu_{_{{\bf r}(\Pi)}}))$\, \ there exists  a vector $ f_{\eta},$
for which $$({\Pi} (g)\eta,\, \  \eta )\, \  =\, \ (\Pi(\sigma_q^{(n)})\Pi (
 g) \Pi ( {(\sigma_q^{(n)})}^*) \Pi(\sigma_q^{(n)})\eta\, \
 ,\Pi(\sigma_q^{(n)}) \eta)\, \ = \, \ $$ $$=\; (\Pi(g_{\sigma} )
 \Pi(\sigma_q^{(n)})\eta\, \ ,\Pi(\sigma_q^{(n)}) \eta)\, \ =\, \
 ({\Pi}_{_{\scriptstyle{Az}}} (g_{\sigma} )\, \  f_{\eta},\, \  f_{\eta}).$$
\medskip

From this follow the classification statements 5.6--
5.7, 5.10.
\medskip

Description of admissible representations of groups
 $ Sp(2\infty)$ \, \         and\, \     $ O(2\infty)$
is essentially based on the structure of the admissible
 representations of \, \  $GL(\infty)$
\, \  and corresponding group of motions.
 As before from the representation
\, \  $\Pi$\, \   of group\, \  $G,$
that coincide here with $ Sp(2\infty)$  \, \
 or\, \    $  O(2\infty),$
we pass to the restriction of\, \  $\Pi$ \, \
 to the subgroup \, \
$GK_n\, \  \subset \, \  G$,\, \  defined
 in \, \  \S 3 \, \
and taking the same part as a group of motions \, \  $G_n^I$ \, \
in a case when $GL(\infty).$
Besides that, this restriction decomposes to the direct integral
 of spherical representations, for which the  explicit realization
 is obtained \, \  (see $ (30), (31), (32),$  (33), (34),
proposition 6.12 and remark 6.13).

Classification concludes (see theorem 6.15) by using theorems 6.14
 and 0.5
as in a case of group $GL(\infty).$

\newpage

\begin{center}
{\Large
\S 1.  Realization Of The Admissible  \newline \hskip 16pt
 Representations Of The
 Symplectic And  \newline  Orthogonal Groups}
\end{center}
\medskip

In this chapter we denote by $G$  one of the groups $Sp(2\infty)$ or
$O(2\infty).$ To define the standard system of generators in  $\, G\, $
for any matrix $\, x\,$ with the elements $\, x_{_{jk}}\, \,
 (j,k=1,2,\ldots)$
we  introduce the matrices
 $\  {\gamma }_{_{\scriptstyle o}}^{(0)}(x) \, \ \,
{\gamma }_{_{\scriptstyle u}}^{(0)}(x) \, \  \in\ G\,$
 by correlations
\medskip

${({\gamma }_{_{\scriptstyle o}}^{(0)}(x)-I)}_{_{jk}}= \left\{
\begin{array}{rl}
x_{_{-jk}}\, ,& \mbox{if}\, j<0\ \, \mbox{and}\, \  k>0\\
0\ \ \ \ \  , & \mbox{otherwise}
\end{array} \right.$
\medskip

${({\gamma }_{_{\scriptstyle u}}^{(0)}(x)-I)}_{_{jk}}= \left\{
\begin{array}{rl}
x_{_{j(-k)}}\, ,& \mbox{if}\, j>0\ \, \mbox{and}\, \  k<0\\
0\ \ \ \ \  , & \mbox{otherwise}.
\end{array} \right.$
\medskip

 Matrices of the form  $g_{_{\scriptstyle 0}}  \, \ =\hskip
4pt \pmatrix{(g^{-1})^{\prime}&0\cr 0&g} \, \ , {\gamma
}_{_{\scriptstyle o}}^{(0)}(x) \, \ , {\gamma }_{_{\scriptstyle
u}}^{(0)}(x) \, \ ,$ where the block structure corresponds to the
decomposition of $H$ into the orthogonal sum of the subspaces $H_-$
and $H_+\,$, generated by basis elements with the
negative and positive numbers respectively, and
$(g^{\prime})_{_{i,k}}=g_{_{-k,-i}}$, are the
system of generators for $G$. As has been above, we identify the
vectors $\,f\,=\,{\sum \limits_i}f_i\,
e_i\, \in H $ with the rows\,$ (\cdots,\, \
f_{-n}\,,\cdots\,, f_{-2}\, ,\,
f_{-1}\,,\,f_1\,,\,f_2\,, \cdots \,)\,$. \medskip

If we denote by $\, t\,$ an ordinary transpose, then the following
 correlations are true:
\medskip

 ${\gamma }_{_{\scriptstyle o}}^{(0)}(x) \, \ =\, \
{\gamma }_{_{\scriptstyle o}}^{(0)}(x^t ) \, \ ,\, \
{\gamma }_{_{\scriptstyle u}}^{(0)}(x) \, \ =\, \
{\gamma }_{_{\scriptstyle u}}^{(0)}(x^t ) \, \ ,\, \  $

  where \, \  $ (x^t)_{kj}\, \  =\, \  x_{(-j)(-k)}$\ \
for \, \  $Sp(2\infty)  \, \ ; $
 \medskip

${\gamma }_{_{\scriptstyle o}}^{(0)}(x) \, \ =\, \
{\gamma }_{_{\scriptstyle o}}^{(0)}(-x^t ) \, \ ,\, \
{\gamma }_{_{\scriptstyle u}}^{(0)}(x) \, \ =\, \
{\gamma }_{_{\scriptstyle u}}^{(0)}(-x^t ) \, \ ,\, \  $

   when \, \ $ G \, \  =  \, \  O(2\infty).$
\medskip

Let $\Lambda _m$ consists of all the complex matrices  $\lambda$ of the form
\newline$\pmatrix{{\lambda}_{11}&{\lambda}_{12}&\ldots&{\lambda}_{1n}&
 \ldots\cr
           {\lambda}_{21}&{\lambda}_{22}&\ldots&{\lambda}_{2n}&\ldots\cr
           \vdots&\vdots&\ddots&\vdots&\ddots& \cr
            {\lambda}_{m1}&{\lambda}_{m2}&\ldots&{\lambda}_{mn}&\ldots\cr}$.
\, \  $\Lambda _m(k) $ \, \   \, \  is a set of columns
\, \  ${\vec{\lambda}}_k\, \  =\, \
\pmatrix{{\lambda}_{1k}\cr {\lambda}_{2k}\cr \vdots\cr {\lambda}_{mk}\cr}\, ,\;
 \, \  A  \, \   \, \  $is a selfadjoint operator on
\, \  $\Lambda _m(k) $ \, \ ,
${\rho}_{_{\scriptstyle k}}({\vec{\lambda}}_k) \, \  = \, \
{ \frac{1}{{\displaystyle {\pi}}^m}} \, \  \exp[-({\vec{\lambda}}_k)^*
 {\vec{\lambda}}_k]\, \ ,
 {\kappa}_{_{\scriptstyle k}}^{^{A}}({\vec{\lambda}}_k) \, \  =
 \, \
\newline \exp[-{i} ({\vec{\lambda}}_k)^* A {\vec{\lambda}}_k]\, \ ,
 \, \
 {\nu}_{_{\scriptstyle m}}^{^{\scriptstyle(k)}}\, \  $ is a measure on
\, \  ${\Lambda}_m(k) $ \, \  with a density   \, \
 ${\rho}_{_{\scriptstyle k}}$ \, \
as regards the Lebesgue measure \newline $ \, \  d{\vec{\lambda}}_k \;,\;
{\nu}_{_{\scriptstyle m}}\, \  =\, \  \prod\limits_{k=1}^{\infty}
{\nu}_{_{\scriptstyle m}}^{^{\scriptstyle(k)}}\,$.
\medskip

We denote by \, \  $L_m^{^A}$ \, \  a Hilbert space,
that is an infinite tensor product \, \
 $\bigotimes\limits_{k=1}^{\infty}
L^2(\Lambda _m(k), d{\vec{\lambda}}_k)$  \, \  with the stabilization
  defined
by the sequence  \, \  ${\eta}_{_{\scriptstyle k}}^{^{A}}\, \
({\vec{\lambda}}_k) \, \  = \, \
{\kappa}_{_{\scriptstyle k}}^{^{A}}({\vec{\lambda}}_k)\, \ \cdot \;
[  {\rho}_{_{\scriptstyle k}}({\vec{\lambda}}_k)   ]^{^{1\over 2}}\,.$
\, \  Namely, \, \  $L_m^{^A}$ \, \  is generated by
 the vectors
of the form \, \  $ f_1 \otimes f_2 \otimes \ldots \otimes f_p \otimes
{\eta}_{_{\scriptstyle k}}^{^{A}} \otimes   {\eta}_{_{\scriptstyle k}}^{^{A}}
 \otimes \ldots  \, \
(p\, \  \in \, \  N \, \ ).$ For $ \, \  f^{(l)}  \, \  = \
   f_1^{(l)}\otimes f_2^{(l)} \otimes \ldots \otimes f_p^{(l)} \otimes
{\eta}_{_{\scriptstyle k}}^{^{A}} \otimes   {\eta}_{_{\scriptstyle k}}^{^{A}}
 \otimes \ldots  \, \  \, \
(l=1, \, \  2)$ the scalar product $\, \  (f^{(1)}\cdot
 f^{(2)})\, \  $ in\, \  $L_m^{^A}$  \, \
is evaluated from formula $$(f^{(1)}\cdot f^{(2)})\, \  = \, \
\prod\limits_{p=1}^{\infty} \int\limits_{\Lambda _m(p)} f_p^{(1)}
 ({\vec{\lambda}}_p ) {\bar f}_p^{(2)}({\vec{\lambda}}_p )
 d{\vec{\lambda}}_pd{\vec{\lambda}}_p\,.$$
\medskip

In$\, \  L^2({\Lambda}_m(k),\, d{\vec{\lambda}}_k) \, \ $
let's define the operator of the Fourier transformation\, \
 $F_k^{^{A}} $:
$$( F_k^{^{A}} f_k) ({\vec{\lambda}}_k)\, \  = \, \
{{\det(\sqrt{\mathstrut{1+4A^2}}\, )}\over{{(2\displaystyle\pi)}^m}}
 \int\limits_{\Lambda _m(k)}\exp\{ i\,\Re\,
 [2\,{({\sqrt{\mathstrut {1+4A^2}} } \, \
 {\vec{\lambda}}_k)}^t \, \  {\hat {\vec{\lambda}}}_k]\}
 \cdot f_k({\hat {\vec{\lambda}}}_k)
d{\hat {\vec{\lambda}}}_k \,.$$
\medskip

The next statement we can obtain using the ordinary calculations.

{\bf Lemma 1.1} {\it If \, \  $A\, \  =\, \  A^* \, \  $
 and
 $A\, \  =\, \ \pm\, A^t \, \ , $
then $\, \  F_k^{^{A}}{\eta}_{_{\scriptstyle k}}^{^{A}} \,=
\, {\eta}_{_{\scriptstyle k}}^{^{A}} \,.$}

\vskip 8pt

From this lemma follows that the operators in the next statement are defined
correctly.

\vskip 8pt

{\bf Proposition 1.2}\, \  {\it Let in\, \   $L_m^{^A}$ \, \
 the action of operators $ \, \  {\Pi}_{_{A}}(g)\, \  (g\, \
\in\, \   Sp(2\infty))$\, \newline  is defined  by formulas :

$$ ({\Pi}_{_{A}}(g_{_{\scriptstyle 0}})\xi )(\lambda) \, \  =\, \
{\vert \det\, g \vert}^m \xi(\lambda g)\, \ ,\, \ $$
\begin{equation}
 ({\Pi}_{_{A}}({\gamma}_{_{\scriptstyle u}}^{(0)}(x) )\xi )(\lambda)\ =
\end{equation}
$$\, \  \exp\{i\,\Re Tr\, [\,
 {\sqrt{\mathstrut {1+4A^2}} } \, \
\lambda x {\lambda}^t]\}\xi (\lambda)\, ,$$
$({\Pi}_{_{A}}(s^-)\xi )(\lambda) \  =\ ( F^{^{A}}\xi)(\lambda),
\, \  $ where $ \;   F^{^{A}}\  = \  \bigotimes\limits_{k=1}^{\infty}
F_k^{^{A}}$.\, \

If  \, \  $A\, \  =\, \  A^* \, \  $and
 $A\, \  =\, \ -\, A^t \,,\   $then operators
$ \;  {\Pi}_{_{A}}(g)\; $ define the representation of  group
\, \ $  Sp(2\infty)$\,.}
\bigskip

Let's give a similar realization for group $O(2\infty).$

{\bf Proposition 1.3.}$\;$ {\it Let \  $ m=2k $ \   and matrix
$\, \  \lambda \, \ = \, \  \pmatrix{{\lambda}^{(1)}\cr
{\lambda}^{(2)}\cr}\, \ ,$ where \, \  $ {\lambda}^{(1)}$\
consists of the first \, $k$   rows of $\lambda \;
\in  \Lambda _m,  P= \pmatrix
{0&{-I_k}\cr{I_k}&0\cr} (I_k\,\;\;\mbox{is}\;\; \, k\times k
\,\, $  unit matrix). If $\, \  A^t\, =
\, \, P\, A  \, P^{-1}\, ,\, $ then
operators  $ \, \  {\Pi}_{_{A}}(g)\, \  (g\, \
\in\, \   O(2\infty)\, ,$\, \   that are defined by the correlations :
$$ ({\Pi}_{_{A}}(g_{_{\scriptstyle 0}})\xi )(\lambda) \, \  =\, \
{\vert \det\, g \vert}^m \xi(\lambda g)\, \ ,\, \ $$
\begin{equation}
 ({\Pi}_{_{A}}({\gamma}_{_{\scriptstyle u}}^{(0)}(x) )\xi )(\lambda)\, \ =
\end{equation}
$$\, \  \exp\{i\, \Re\, Tr\,  [\, {\sqrt{\mathstrut {1+4{(A^t)}^2}} } \
P \lambda x {\lambda}^t]\}\xi (\lambda)\,, $$
${\Pi}_{_{A}}(s^+) \, =\,  F^{^{A^t}},$
define the representation of group $O(2\infty).$}
\medskip

For the decomposition of the representation  $\ {\Pi}_{_{A}} \ $
 to the irreducible
components let us consider two groups :
\begin{equation}
O(A, m) = \{ u\in U(m) : u^t = u^*
\, \  and \;\; [ A, u]
= 0\, \}\, ,
\end{equation}
\begin{equation}
Sp(A, m) = \{u\in
U(m) :   u^t = P u^* {P^{-1}}
\, \  and \;\; [ A, u]
= 0 \}.
\end{equation}

Let \  $( \tau (u) \xi)( \lambda)\  =\  \xi(u^{^{-1}} \lambda)\, ,$
 where $\, u\, \in\, U(m)\, ,\, \xi\, \in\,  L_m^{^A}.$
\medskip

If \  $G(A,m)\;$  is one of the groups\, $ O(A\, ,\, m)\ $
or \ $ Sp(A\, ,\, m)\, ,\  \rho \, $ is its irreducible representation,
 $\, a_{_{kk}}^{\rho}\, $\,
is matrix element of \, $\rho \, , $\,  then
$  P_k^{\rho}  =  \dim\,\rho
\int\limits_{G(A,m)} \, a_{_{kk}}^{\rho}(u)\tau (u) \, du
 \,\;$ is an orthogonal projection, that acts in $\,  L_m^{^A}.$
\medskip

The following statement was announced at first by G.I. Ol'shansky in
 \cite{Olsh3}.

{\bf Proposition 1.4.} {\it  Let \,
 $ {\wp}$\,
 be some set of the operators, that act on $\,  L_m^{^A},$\
$ {\wp}^\prime$\, is a commutant of \ $ {\wp}$\,.\,
Then the following statements are true:
\begin{center}
{\sc i)} if  $ {\Pi}_{_{A}}$  such, as in
proposition 1.2, then $${[{\Pi}_{_{A}}(Sp(2\infty))]}^\prime\  =\
{({(\tau(O(A\, ,\, m)))}^\prime)}^\prime =
{{\tau(O(A\, ,\, m))}^\prime}^\prime ;$$
{\sc ii)} if  $ {\Pi}_{_{A}}$ such, as in
proposition 1.3, then $${[{\Pi}_{_{A}}(O(2\infty))]}^\prime  =
{{\tau(Sp(A\, ,\, m))}^\prime}^\prime ;$$
{\sc iii)} restriction of $ {\Pi}_{_{A}}$ to
 $ P_k^{\rho}\,  L_m^{^A}$
is irreducible.
\end{center}  }
\vskip 10 pt
\begin{Large}
\begin{center}

\S 2 \ \ \  Some Properties  Of The Representation\newline Of Groups $GL(\infty)$
 \end{center} \end{Large}
\vskip 8pt

Let $\Pi$ be an irreducible representation of  $GL(\infty)$, acting
in a Hilbert space
$H_{_{\Pi}}$ \  with the unit cyclic vector $\xi $
 fixed with respect to the operators \  $\Pi(u)\  (u\, \in \
 U(GL(\infty))\, =\,  U(\infty)).  $

Then by the results of the  paper \cite {N1}
the class of the unitary equivalence of  $\Pi\; $ is defined by
 a spherical function
$${\varphi}_{_{\Pi}} = (\Pi(g)\xi, \xi) =$$
 \begin{equation} {|\det(g)|}^{i\beta} \det[I_{_{{\bf r}(\Pi)}}
\otimes ch \ln |g| - 2 i A \otimes sh \ln |g|],
\end{equation}
where $|g|\ =\ (g^*g)^{1\over2}\ \ ,\
 {\bf{\scriptstyle r}}({\scriptstyle{\Pi}})$ is a natural number,
  $\, \ \beta $
is a real number, $I_{_{{\bf r}(\Pi)}}\; $ is a unit\ \ \ \ \ \
 ${\bf{\scriptstyle r}}({\scriptstyle{\Pi}})\times
 {\bf{\scriptstyle r}}({\scriptstyle{\Pi}})-$ matrix,
 $A$ is selfadjoint (s.a.) ${\bf{\scriptstyle r}}({\scriptstyle{\Pi}})\times
{\bf{\scriptstyle r}}({\scriptstyle{\Pi}})-$ matrix.
\medskip

{\sl The natural number ${\bf{\scriptstyle r}}({\scriptstyle{\Pi}})$
 we call the rang of the representation $\Pi$.}
\medskip

We denote by \ $Z_m$ \ a set$\ m\times m$ of the upper triangular matrices
with  positive elements on the diagonal.
Let subgroup\ $D_m\ \subset\ GL(\infty)$\  consists of the matrices
 of the form
$\pmatrix{z_m&0\cr x_m&I\cr} \,,$ where $z_m\ \in \ Z_m\,,\; x_m\ $ is
 matrix with the  finite number of nonzero elements.  \medskip

Let us give a classification result for the spherical
 factor representations of group
\ $G_n^I\ \subset \ GL(\infty)\,$ that consists of matrices
 of the form\ $\pmatrix{I_n&0\cr
*&*\cr}\,,$ where \ $\ast$\ \  signifies some matrix of the corresponding
 size.

{\bf Theorem 2.1.} (see \cite{N2}) {\it \ \ Let $\Pi$ be irreducible
 spherical factor representation
\ $G_n^I\,$ that acts in\ $H_{_{\Pi}} \,$\
 ${\hat{\alpha}_{_A}}(\lambda,g) = \exp \lbrace \frac {-1}{2}
Tr \lbrack \lambda(gg^{\ast}-1)\lambda^{\ast}-
2i A\lambda(gg^{\ast}-1)\lambda^{\ast} \rbrack \rbrace, $ where $A\ $ is
s.a. $m\times m\ -$ matrix.

Then there exist : s.a. $m\times m\ -$ matrix $A\,,\ n\times m\ -$ matrix
$z$\,, real number $\beta\ $
such, that $\Pi$ unitary equivalent to the restriction of representation
${\Pi}_{_{\scriptstyle{Az}}}$\ \ of group $G_n^I\,$  defined in
\ $L^2(\Lambda_m,\nu_m)$ \ by formulas:
$$({\Pi}_{_{\scriptstyle{Az}}}{\bigg (}\pmatrix{I_n & 0 \cr 0 & g}
 {\bigg )}\eta)(\lambda)= \mid det  g \mid ^{i\beta}
{{\hat{\alpha}_{_A}} } (\lambda,g)\eta(\lambda g), $$
\begin{equation}
({\Pi}_{_{\scriptstyle{Az}}}{\bigg (}\pmatrix{I_n & 0 \cr h & I}
 {\bigg )}\eta)(\lambda)=
\exp \lbrack i \, \  \Re \, \  Tr (z \lambda h ) \rbrack
\eta (\lambda)
\end{equation}
$$( \eta \in L^2 ( \Lambda_m, \nu_m)),$$
to the subspace\ \ $[{\Pi}_{_{\scriptstyle{Az}}}(G_n^I)\xi_{_{0}}]\,\ , $
where $\xi_{_{0}}\ \in\ L^2 ( \Lambda_m, \nu_m) $ and corresponds to the
 function on  $\Lambda_m\,$ that identically equals to the unit.} \vskip
8pt

Using the form of the operators \ ${\Pi}_{_{\scriptstyle{Az}}}(g)\ \
 (g\ \in\ G_n^I)\,$
we can prove the following statement.
\medskip

{\bf Theorem 2.2.} \ \ {\it Let subgroup  $D_{mn} \subset G_n^I\,  (G_0^I \ =
 \, GL(\infty)\, ,\, D_{m0} \  =\  D_m)  $
consists of matrices of the form
 $\pmatrix{I_n & 0 & 0 \cr * & z_m & 0 \cr * & * & I \cr}\ $\
 $ (z_m\ \in \ Z_m)\, ,\ \Pi  \  $  \
is the same as in theorem 2.1. \
 Then $ [\Pi(D_{mn})\xi]\ =\ [\Pi (G_n^I)\xi]\,,$
where\ $\xi\ \ $ is an arbitrary  nonzero \ $  \Pi (U(G_n^I))\,-$  fixed
vector.}

\vskip 16pt
\begin{center}
{\Large
\S 3.   Asymptotic Properties \newline
 Of The Admissible Representations}
\end{center}

\medskip

Let $G\, =\,  Sp(2\infty)$ or
$O(2\infty)\,\,$. The subgroups\ $G_d\, ,\ \ \Gamma_{_O}\, ,\
 \Gamma_{_U}\ \subset \ G$\ \
consist of matrices of the form \ $g_{_{\scriptstyle 0}}  \hskip
4pt=\, \  \pmatrix{(g^{-1})^{\prime}&0\cr 0&g} \, \ , {\gamma
}_{_{\scriptstyle o}}^{(0)}(x) \, \ , {\gamma }_{_{\scriptstyle
u}}^{(0)}(x) \, \ $ respectively, where $\, x\, $is a matrix with the
elements $\, x_{_{jk}}\, \, (j,k=1,2,\ldots)\,,$ and the block structure is
defined by the decomposition of $H$ into the orthogonal sum of the subspaces
$H_-$ and $H_+\ $ (see \S 1).  \medskip

If ${(x^t)}_{_{kj}}\ =\ x_{_{jk}}\,$, then the correlations are true
\begin{eqnarray}
{\gamma }_{_{\scriptstyle o}}^{(0)}(x) \, =
 \, {\gamma }_{_{\scriptstyle o}}^{(0)}(x^t ) \, ,
 \ {\gamma }_{_{\scriptstyle u}}^{(0)}(x) =
 {\gamma }_{_{\scriptstyle u}}^{(0)}(x^t ) \ \mbox{for}  \  Sp(2\infty);
 \ &\nonumber \\
{\gamma }_{_{\scriptstyle o}}^{(0)}(x) \, =
 \, {\gamma }_{_{\scriptstyle o}}^{(0)}(-x^t ) \, ,
 \, {\gamma }_{_{\scriptstyle u}}^{(0)}(x) \, =
 \, {\gamma }_{_{\scriptstyle u}}^{(0)}(-x^t ) \, \mbox{for}  \, G \, =
 \, O(2\infty)&.
\end{eqnarray}
\medskip

Let $\, b\,$ be a matrix with the elements $\, b_{_{jk}}\, ,$ where
 $\, j,k=1,2,\ldots \,.$
As in \S 1 let's define the matrices
 $ {\gamma }_{_{\scriptstyle o}}^{(n)}(b) \, $ and
 ${\gamma }_{_{\scriptstyle u}}^{(n)}(b) \ \in\ G\,$
by the correlations
\newline
\medskip

${({\gamma }_{_{\scriptstyle o}}^{(n)}(b)-I)}_{_{jk}}= \left\{
\begin{array}{rl}
b_{_{(-j-n)(k-n)}}\, ,\;& \mbox{if}\, j<-n\ \, \mbox{and}\, \  k>n\\
0\ \ \ \ \  ,& \mbox{otherwise}\,
\end{array} \right.$
\medskip

${({\gamma }_{_{\scriptstyle u}}^{(n)}(b)-I)}_{_{jk}}= \left\{
\begin{array}{rl}
b_{_{(j-n)(-k-n)}}\, ,\;& \mbox{if}\, j>n\ \, \mbox{and}\, \  k<-n\\
0\ \ \ \ \  ,& \mbox{otherwise}.
\end{array} \right.$
\medskip

Let $Y(\infty, n)\, $ be a set of all the local nonzero
 $(\infty\times n\, )$-matrices
with the infinite number of rows and $n\, $ columns.
 If $\, x=[x_{_{jk}}]\, ,\,
y=[y_{_{jk}}]\ (1\leq j <\infty\, ,\, 1\leq k\leq n)\, ,$ then  we define
matrix ${\theta}^{(n)}(x,y)\, \in G\,$ by the correlation
\newline
\medskip

${({\theta}^{(n)}(x,y)-I)}_{_{im}}\, =\left\{
\begin{array}{rl}
x_{_{(i-n)m}}\, ,\;& \mbox{if}\, i>n\ \, \mbox{and}\, \  1\leq m\leq n \\
-(x^t)_{_{(-i)(-m-n)}}\, ,\;& \mbox{if}\,\ -n\leq i\leq -1\ \, \mbox{and}\,
 \   m<n \\
y_{_{(-i-n)m}}\, ,\;& \mbox{if}\,\ i<-n\ \, \mbox{and}\, \ 1\leq  m\leq n \\
{y^{\sharp}}_{_{(-i)(m-n)}}\, ,\;& \mbox{if}\,\ -n\leq i\leq -1\ \, \mbox{and}\,
 \ m>n \\
0\ \ \ \ \  ,& \mbox{otherwise}\, ,
\end{array} \right.$
\medskip

where $\, {y^{\sharp}}\, =\, y^t \ \ for \ \ \ G \, =\,
Sp(2\infty) \ and \ y^{\sharp}= -y^t \ \
for \ \ \  \, G  =   O(2\infty)$.
\medskip

At last we denote by $\, {\delta}^{(n)}(a)\ $
an element from $\, G\, ,$  that defined by $\, n\times n\, -$matrix
 $\, a=[a_{_{jk}}]\ (1\leq j,k \leq n)\, $
with respect to the correlation
\medskip

$\,{ ({\delta}^{(n)}(a)-I)}_{_{im}}=\left\{
\begin{array}{rl}
a_{_{(-i)m}}\, ,\;& \mbox{if}\, -n\leq i\leq -1\ \, \mbox{and}\, \
 1\leq m\leq n \\
0\ \ \ \ \  ,& \mbox{otherwise}.
\end{array} \right.$
\medskip

Let
\medskip

$g_n  =
\pmatrix{{(g^{-1})}^{\prime} &0&0&0 \cr 0& I_n &0&0\cr
 0&0& I_n &0\cr 0&0&0&g\cr}\, ,$
where $\, {(g^{\prime})}_{_{ik}}=g_{_{(-k)(-i)}}\, \, (k,i<0)$.
\medskip

Let\ $G_d(n,\infty)\ (G_d(0,\infty)\, =  \, G_d),\ \Gamma_{_O}^{(n)}
\ (n\geq 1),\ \Gamma_{_U}^{(n)}\ (n\geq 1),\ \,\
\Delta^{(n)}\ (n\geq 1)\ \,$ are the subgroups of  $G$,\ which consist
 of the matrices of the form $g_n\, ,\,
{\gamma }_{_{\scriptstyle o}}^{(n)}(b)  ,\
 {\gamma }_{_{\scriptstyle u}}^{(n)}(b)  ,\,
 {\delta}^{(n)}(a) $\ \  respectively.
\  Moreover, when $G\,  =  \, Sp(2\infty) \, ,\ a $ is symmetric matrix,
 and when
$ O(2\infty)\ -$  is antisymmetric one. A set of all the elements
from $G\, $  of the form $\,  {\theta}^{(n)}(x,y)\, $ we denote by $\,
\Theta^{(n)}\ (n\geq 1) \,.$ Observe that $\,  \Theta^{(n)}\ $ is not a
group.  \medskip

The following correlations are true :
\begin{eqnarray}
&&{\gamma }_{_{\scriptstyle o}}^{(n)}(b) {\theta}^{(n)}(x,y)\
 {\gamma }_{_{\scriptstyle o}}^{(n)}(-b)=\nonumber \\
&&{\delta}^{(n)}(-{(bx)}^{\sharp}x)\ {\theta}^{(n)}(0,bx){\theta}^{(n)}(x,y)
 \,\nonumber \\
&&{\gamma }_{_{\scriptstyle u}}^{(n)}(b){\theta}^{(n)}(x,y)\
 {\gamma }_{_{\scriptstyle u}}^{(n)}(-b) =\nonumber \\
&&{\delta}^{(n)}({(by)}^ty)\ {\theta}^{(n)}(by,0)\ {\theta}^{(n)}(x,y)\, \\
&&{\theta}^{(n)}(x_{_1},y_{_1})\ {\theta}^{(n)}(x_{_2},y_{_2})=  \nonumber \\
&&{\delta}^{(n)}(x_{_2}^t y_{_1} - y_{_2}^{\sharp} x_{_1} - x_{_1}^t y_{_2} +
  y_{_1}^{\sharp} x_{_2} ) {\theta}^{(n)}(x_{_2},y_{_2}){\theta}^{(n)}
 (x_{_1},y_{_1})\, ,\nonumber \\
&&g_n\ {\theta}^{(n)}(x,y)\ g_n^{-1} =
 {\theta}^{(n)}(gx,{(g^t)}^{-1}y).  \nonumber
\end{eqnarray}

 Let  subgroup  $K_n$ \ be generated by $\Theta^{(n)}\ $ and
 \ $\Delta^{(n)}$.
From the correlations \ (8) \ it follows, that the elements of the group \
 $G(n,\infty)$, \ that is formed by
\ $G_d(n,\infty)\,,\  \Gamma_{_O}^{(n)},
\ \Gamma_{_U}^{(n)}\, $ naturally act   by the automorphisms on\ $K_n$.
 We denote by
\ $GK_n$ \ a subgroup of  $G$ generated by \ $G(n,\infty)$ \
 and\  $K_n$.

Let's  introduce a useful property of the infinite-dimensional groups.
\medskip

Suppose, \ that group\ $G$ \ is an inductive limit of the finite-dimensional
matrix groups\ $G(n)$ \ $(G(n)\subset   G(n+1)\;  (n=1,\ 2,\ 3,\
 \ldots ))$.
\medskip

{\bf Definition 3.1.}{\it \ \ \ Group \ $G$ \ is called {\sl an asymptotic
 Abelian one
\ \ (a.a.)}, if there exists a sequence \ $u_n\ (n\in N)$ \ of the elements
of  subgroup $ \ U(G)\,$   that coincide with a set of the
 unitary  elements of the group
 \ $G\,$  \ with the properties :
\newline

{\sl i)} for any  $\l,\ n\ \in\ N\ $   and $\ g\in G(l)$
 there exists  $k(g,n)\ \in\ N$\ such that  for $k\geq k(g,n)
\ \ u_k\ g\ u_k^*\  \in\ \{h \in G : \ hp=ph \  $ for   all $\ p \in G(n)\}\ =
 \ G^{\prime}(n)$;
\newline

{\sl ii)} for any  $n\in N$  there exists $l(n)\ \in\ N, \ $
 for which
$u_k \,  u_l^* \ \in\ G^{\prime}(n)\,$ \  when $k>l\geq l(n)$}.
\medskip

{\bf Proposition 3.2.}{\it \ \ \ Set \ $G$ is Á.Á. matrix group,
 that is an inductive limit of the finite-dimensional groups
$G(n)$ \ $(n\in N)\, ,$ \ $\Pi $ is a factor--representation  of $G\,
 ,$   that acts in a Hilbert
space $H_{_\Pi}$ with a cyclic vector $\xi$. Moreover, set
\newline
$\bigcup\limits_{n=1}^{\infty}
\{\eta \in H_{_\Pi}\ :\ \Pi (u) \eta\ =\ \eta$  for all
$u\in U(G)\ \cap \ \{h \in G\ :\ h\ g\ =\ g\ h\ \ $ for all \ $g \in G(n)\}\}$
dense in \ $H_{_\Pi}$\ ({\sl This property, as was shown by G.I. Ol'shansky,
 is one
of the equivalent definitions of the admissible representation in a case of
the large class of the infinite-dimensional groups, that includes}
 $GL(\infty)\, ,\ Sp(2\infty)\, ,\;
O(2\infty)\ )$.
\medskip

ôhen for any unit vectors $\ \eta_1,\, \eta_2\ \in \ H_{_\Pi}$

$$\varphi_{{\eta}_1}(g)\, =\,
\lim_{l\to \infty}(\Pi (u_l\ g\ u_l^*)\ \eta_1\, ,\, \eta_1)\ =
\ \lim_{l\to\infty}(\Pi (u_l\ g\ u_l^*)\ \eta_2\, ,\, \eta_2)\ =\
\varphi_{{\eta}_2}(g)\,.$$

Moreover, the limits do not depend on a choice of the sequence
 $u_l\ (l\in N)\ $
from the definition 3.1,\ and\ $\varphi_{{\eta}_1}\, =\,
 \varphi_{{\eta}_1}\ $
is an indecomposable  spherical function on\  $G$\ .}
\medskip

{\bf Proof}. For any $\epsilon >0$ \ there exists $n\in N$ and
unit vectors $\xi( \epsilon )\,\ \eta_i (\epsilon) \in \ H_{_\Pi}\ \
 (i=\ 1,\ 2)\ \ $
with the properties:
$\Pi (u)\xi( \epsilon )= \xi( \epsilon )\, ,\;
\Pi (u) \eta_i (\epsilon) = \eta_i (\epsilon)
\ (i=1,\ 2)$  for all  $u \in\ U(G)\cap  G_n^{\prime}\ ;$
\begin{equation}
||\xi\ -\ \xi( \epsilon )||< \epsilon\ ;\ ||\eta_i \ -\ \eta_i
 (\epsilon)||<\epsilon \ .
\end{equation}

By the definition 3.1 we choose $l\in N $  such, that\  $u_k \,
u_l^*\ \in\ G_n^{\prime}$ \ for all \ $ k>l $. Then, using  (9)  we get :
$\newline
2 \epsilon >  |(\Pi (u_k\ g\ u_k^*)\ \eta_i\,\ \eta_i)\ -
(\Pi (u_k\ g\ u_k^*)\ \eta_i({\epsilon})\,\ \eta_i({\epsilon}))|\ =
\newline
\medskip
|(\Pi (u_k\ g\ u_k^*)\ \eta_i\,\ \eta_i)\ -
(\Pi (u_l\ g\ u_l^*)\ \eta_i({\epsilon})\ ,\ \eta_i({\epsilon}))|\ .$
\medskip

Hence and from  (9)  it follows that
$\newline
|(\Pi (u_k\ g\ u_k^*)\ \eta_i\,\ \eta_i)\ -\ (\Pi (u_l\ g\ u_l^*)\ \eta_i\,
 \ \eta_i) |<4\epsilon $
for all $k>l$.

Since   $\epsilon $  is an arbitrary,
 then $ (\Pi (u_k\ g\ u_k^*)\ \eta_i\,\ \eta_i)$ is
a fundamental sequence.  Therefore the limits of sequences
  $ (\Pi (u_k\ g\ u_k^*)\ \eta_i\,\ \eta_i) \ (k\in N) $
 exist for every
 $i=1,\ 2$.

Moreover, from the correlation
$$\lim_{l\to \infty}\Pi (u_l\ u\ u_l^*)\ \eta_1\, =\, \eta_1
\ \ \forall \ u\in U(G)$$
we get, that
$$\varphi_{{\eta}_1}(ugv)\, =\, \varphi_{{\eta}_1}(g)\,  \
 \forall \ u,\ v \in U(G)  \ .$$
Therefore  $\varphi_{{\eta_1}}$ and $\varphi_{{\eta}_2}$
 are the spherical
functions on $G$.

Let's prove the coincidence  of  $\varphi_{{\eta}_1}$  and
 $\varphi_{{\eta}_2}$.

As  $\xi $ is a cyclic vector, then for any  $\delta >0 $   there exist
the collections  $ {\{ g_{_{{\displaystyle ik}_{\scriptstyle i}}} \} }_{k_i =
 1}^{N_i}\ \ (i=1,\ 2)$ of the elements from
\ $G(n)$   and  numbers
 ${\{ c_{_{{\displaystyle ik}_{\scriptstyle i}}} \} }_{k_i =1}^{N_i}\ \
 (i=1,\ 2)$ from
 ${\bf C} $ with the properties
\begin{equation}
||\sum_{k_i=1}^{N_i}\   c_{_{{\displaystyle ik}_{\scriptstyle i}}}\
 \Pi ( g_{_{{\displaystyle ik}_{\scriptstyle i}}}) \xi \
-\ \eta_i\ ||\ <\   \delta \ \ \ \ (i=1,\ 2)
\end{equation}

Since  $\Pi $  is factor-representation, then, using  (10), we get
$$ 2\delta > \ |\lim_{l\to \infty} (\Pi (u_l\ g\ u_l^*)\ \eta_i\,\ \eta_i)
\ -\ $$ \ $$\lim_{l\to \infty}\ (\Pi (u_l\ g\ u_l^*)\
\sum_{k_i=1}^{N_i}\,  c_{_{{\displaystyle ik}_{\scriptstyle i}}}\
 \Pi ( g_{_{{\displaystyle ik}_{\scriptstyle i}}}
) \xi \,\ \sum_{k_i=1}^{N_i}\,  c_{_{{\displaystyle ik}_{\scriptstyle i}}}\
 \Pi ( g_{_{{\displaystyle ik}_{\scriptstyle i}}}
) \xi \ )|\ = $$$$\ {\bigg |}\lim_{l\to \infty} (\Pi (u_l\ g\ u_l^*)\ \eta_i
 \,\ \eta_i)
\ -\ \lim_{l\to \infty} (\Pi (u_l\ g\ u_l^*)\ \xi\,\ \xi)\
 {{\bigg |\bigg |} \sum_{k_i=1}^{N_i}\,
 c_{_{{\displaystyle ik}_{\scriptstyle i}}}\
 \Pi ( g_{_{{\displaystyle ik}_{\scriptstyle i}}}
) \xi \ {\bigg |\bigg |}}^2\bigg | \ .$$

Since  $\ \delta \ $ is arbitrary, then from this and from  (10)
 follows, that
$\,\varphi_{{\eta}_1}(g)\, =\ \,\varphi_{{\eta}_2}(g)\,
\ \ \forall \ g\in G \ .$

The proof of the indecomposability of $ \varphi_{{\eta}_1} $
 is based on the correlation
  $$\ \lim_{l\to \infty}\ \varphi_{{\eta}_1}(u_l g\ u_l^*\ h)\, =\ \,
 \varphi_{{\eta}_1}(g)\ \varphi_{\eta_1}(h)\,
\ \ \forall \ g,\ h\in G \ .$$
\medskip

{\bf Definition 3.3.  } {\it A spherical function on  a.a. group
$G$  defined in accordance with the proposition 3.2 by  the
admissible representation $\Pi$ is called an asymptotic spherical
function (a.s.f.)  and it is denoted by  $\varphi_{_{\Pi}}^{(a)}$.}
\medskip

{\bf Theorem 3.4.} {\it Let $G$ be an a.a. matrix group, that is an
inductive limit of finite-dimensional groups   $G(n) (n\in N)\,,\Pi$ is
the same, as in the proposition 3.2, $\xi$ is a unit cyclic vector for $\Pi
\,$ fixed with respect to  $\Pi(u_l)\:(l\in \ N),$ where  $u_l$ is a
sequence from $U(G)\,$ that guarantees a.a. of $G.$

ôhen subspace  $H_{_U}\ =\ \{\eta \in H_{_{\Pi}}\ : \ \Pi(u_l)\eta
\ =\ \eta\ $ for all  $l\in \ N \}$  is one-dimensional.}

{\bf Proof} The proof is completely similar to the proof of the
 {\sl  Multiplicativity Theorem} from  \cite{N1}.
\medskip

{\bf Definition 3.5.}  {\it We define the  {\sl rang}
 ${\bf r}(\Pi )$ of the admissible factor--representation
 $\Pi $  of group  $\ G,$  that coincides to the one of the
 following groups  :
$\ GL(\infty)\, ,\; G_n^I\, ,\; Sp(2\infty)\, ,\; O(2\infty)\, ,\;\newline
GK_n\,$  as a rang of the representation  ${\Pi}^{(a)}\,$ that defined by
the indecomposable Á.s.f.  $\varphi_{_{\Pi}}^{(a)}\,$  of group
 $GL(\infty)\ (G\ =\ GL(\infty))\, ;\;
GL(n, \infty) \,, $ where $\ G\ =\  G_n^I\, ;\; \ G_d\ (G\ =\ Sp(2\infty)\ $\
 or \ $
O(2\infty)\ )\, ;\; \  G_d\cap G(n, \infty)\  (G\ =\ GK_n)\ $}(see \S 2).

\bigskip

\begin{center}
{\Large
\S 4 .   Properties Of The Admissible Representations   \newline
                    Of Group$\ GL(\infty)\ $ \ And The Group Of Motions }

\end{center}
\medskip

We may take the group of motions to be the subgroup of
 $G_n^I\in \ GL(\infty)\ $
(see \S 2)
\medskip

{\bf Theorem 4.1. } {\it Let $\Pi$ be an admissible representation of $G\, ,$
 where  $G$ coincides with the one of the  groups $GL(\infty)\, ,\; G_n^I \,
,\;Sp(2\infty)\,$,\ $O(2\infty)\, ,\; GK_n\ ;\ {\xi}_1 \,,\; {\xi}_2
\;\mbox{are} \;\Pi (U(G(n,\infty ))\ - $fixed vectors from  $H_{_{\Pi}}$.
ôhen $(\Pi(g){\xi}_1\,\,{\xi}_1)\, =\, (\Pi(g){\xi}_2\,\,{\xi}_2)\,$ for all
$g\,\in\,G(n,\infty)\,$.}

{\bf To prove}  this theorem it suffices to  notice, that the groups
 from  the condition are the asymptotic Abelian
ones, the set\newline    $\bigcup\limits_n\,\{\eta \in H_{_{\Pi}}\, :
 \, \Pi (U(G_n^{\prime}))\eta\,=\,\eta \}$
is dense in $ H_{_{\Pi}}\  $ and to use the proposition 3.2.

\vskip 6pt

In $GL(\infty)\ $ we consider a subset  $Y_p \,$ that consists of
 matrices of the form
$$\pmatrix{y_{11}&y_{12}&\ldots&y_{1p} &{\delta}_1 &0&0 &\ldots &0&\ldots&
 \ldots \cr
           *&*&\ldots  &   *&  * & {\delta}_2& 0&\ldots &0 &\ldots &\ldots\cr
         \vdots &\vdots&\vdots &\vdots& \vdots &\vdots&\ddots&\vdots&
 \vdots&\vdots&\vdots \cr
*&*&\ldots& *&*& *& *&{\delta}_s&0&\ldots&\ldots\cr
\vdots&\vdots& \vdots&   \vdots&\vdots&\vdots&\vdots&\vdots&\vdots&
 \vdots&\vdots\cr} \;\; \, ,$$
where $p\geq 1\, ,\ {\delta}_s\, >\,0 \ \   \forall\ s\in N.$
\medskip

{\bf Lemma 4.2.} {\it Let $\Pi $ be an admissible factor--representation of
group $ GL(\infty)\,,\ \xi \ $ is its cyclic vector fixed with respect
 to the operators  $\Pi(U(G(p,\infty))) \ .$
 ôhen  ${\big [}\Pi(Y_p)\xi {\big
]}\, = \,{\big [} \Pi(GL(\infty)) \xi {\big ] }.$}

{\bf To prove} the lemma it suffices to  notice, that
 $GL(\infty)\, =\, Y_p\, U(G(p,\infty))\ .$
\medskip

{\bf Lemma 4.3.} {\it Let $\Pi\, ,\, p\, ,\ \xi\ $ are the same as
in   lemma 4.2,  $\ K_{2p}\ $is a subgroup of $GL(\infty)\,$
 that consists of matrices of the form
$\pmatrix{g_{11}&0\cr *&I\cr} \, , $ where  $\ g_{11}\ $
 is an arbitrary  $\ {(2p)}\times {(2p)} \, - $
matrix, $\ p\,\geq\, {\bf r}(\Pi )\,.$ ôhen $H_{_{\Pi}}\,=
 \,{\big [}\Pi(K_{2p}) \xi{\big ] }\ .\ $}

\medskip

{\bf Proof.} It's easy to check it, that the closing
 of the set  $\ K_{2p}GL(p,\infty )\ $
includes $\, Y_p\,.$ Next, by the statements 3.2 and 4.1
 a restriction of $\Pi $ to
 $\, GL(p,\infty )\, , $ that acts in a Hilbert space
 $[\Pi ( GL(p,\infty )\, )\xi ]\,, $
is a factor--representation. According to the  theorem 2.2
 (for $ n=0 $)  $\,  [\Pi ( GL(p,\infty ) )\xi ]\, = \,
[\Pi ( GL(p,\infty )\cap K_{2p} \, )\xi ]\,$.
 From this and from lemma 3.2  we get
 $\, [{\Pi} ( GL(\infty ) )\xi ] \, = \,
 [ \Pi (Y_p ) \xi \, ] \, = \, [\Pi (K_{2p} GL(p,\infty )) \xi ]\, =\,
[\Pi (K_{2p}) \xi ] \,.$ Lemma 4.3 is proved.

\vskip 6pt

 The following analog of the statement above might be proved
by the same way  for group
\, $G_n^I \subset GL(\infty ) \,$.
\medskip

{\bf\ Lemma 4.4. } \ \ {\it Let $\Pi \,$ be an admissible
 factor--representation
of group $ G_n^I \ \,, \xi\, $ is a cyclic vector fixed
 with respect to
operators $\Pi (u)  (  u\in U(G(p+n,\infty) )\, ,\ K_{2p}^{(n)}\, $
 is a subgroup of
  $\; G_n^I \, ,$ that consists of the matrices of the form
$ \pmatrix{I_n & 0 & 0 \cr
            * & g_{_{\scriptstyle 22}} & 0 \cr
            * & * & I \cr}  \ \, ,$
where  $g_{_{\scriptstyle 22}} \, $ is  an arbitrary
 $\, {(2p)}\times {(2p)} \, - $ matrix.

If \, $p\geq {\bf r} (\Pi ) $ \, (see definition 3.5 ),
 then\, $[\Pi (  K_{2p}^{(n)}\, ) \xi ]\,=
\, [\Pi (G_n^I) \xi]\,.$}

\vskip 12pt
\begin{center}
{\Large
\S 5  .      A Description Of The Admissible

\hskip 16pt  Representations  Of The Groups

\hskip 16pt $\ GL(\infty)\ $ \ And \ $G_n^I $}

\end{center}

\vskip 4pt

Let $\Pi \, $ be an admissible factor--representation of  $G\, ,$
that in  first two statements of this chapter may be one of the groups:
 $\ GL(\infty)\, ,\ G_n^I \, ,\ Sp(2\infty)\, ,\ O(2\infty)\, .$
\medskip

Isometry  $\sigma_q^{(n)}  $  \, ,\ that acts in $H$ by the formula:
$$\sigma_q^{(n)}(e_i) = e_i\;\;\mbox{{\rm for}}\;\; i\leq n
\;\;\mbox{{\rm and}}\;\;
   \sigma_q^{(n)}(e_i) = e_{i+q}
\;\;\mbox{{\rm for}}\;\; i> n $$\ $(n=0\,\;$ corresponds to $\; GL(\infty))
,$ is a strong limit of the elements of unitary subgroup $U(G_n^I) \, .$
 In a case
of groups \, $ Sp(2\infty)\, ,\ O(2\infty)\, $  this property belongs to
  isometry
 \, $ \sigma_q \,  $  defined by the correlation \,
 $ \sigma_q(e_i)\, = \,
e_{i+(sign\, i)q}\, .$

\medskip

We denote by ${\Pi }_{_U} $ a restriction of $\Pi$ to $U(G)  .$
 By the results
of paper \cite{Olsh1} ${\Pi }_{_U} $ is a continuous, if
 $U(G)\, \subset \, B(H) $\ and\  $U(H_{_{ \Pi}})\, \subset \,
 B(H_{_{ \Pi}})$ are the topological groups
with respect to the strong operator topology.\
 Therefore, \,${\Pi }_{_U} $
extends by the continuity to the representation of
 semigroup ${\breve U}(G)\, ,$
that includes all the isometries, which are the limiting points
with respect to the weak operator topology  on $U(G)  .$
Obviously, $  \sigma_q^{(n)}\, \, ,\,\,  \sigma_q $ \
 belong to \ ${\breve U}(G)\, .$
\medskip

{\bf  Theorem 5.1. } \ \ {\it If $E_l \, $is an orthoprojection to the
subspace \ $\Pi (\sigma)\, H_{_{ \Pi}} \, , $ where $\sigma $ equals to
 $\sigma_q^{(n)}$ or $ \sigma_q \, ,$
then \ $  \Pi (g)\,E_l \, - \, E_l \, \Pi (g)\, =
 \, [\, \Pi (g)\, ,E_l \, ]\, =
\, 0 \, $ for all $ g \in G(l,\infty )\, .$}
\medskip

{\bf Proof.} Since $\Pi $ is extended by the continuity to \  $\sigma $\
\, ,\, then for arbitrary\ \,  $g\in G(l,\infty) \,\,\,\,
 {\sigma}^*\, g\, {\sigma }\, \in \, G\,$ and
\, $ \Pi (  {\sigma}^* )\, \Pi (g)\, \Pi ({\sigma })\, = \,
 \Pi (\, {\sigma}^*\, g\, {\sigma }\, )\, .$
\medskip

Therefore, \, $ \Pi (  {\sigma}^* )\, E_l\, \Pi (g)\, E_l\,
 \Pi ({\sigma })\,
$is a unitary operator in  \, $H_{_{ \Pi}} \, .$ Hence, accounting that
$\Pi (\sigma ) )\, $ is an isometry, we get, that  $\, E_l\,
 \Pi (g)\, E_l\, \, $
is a unitary operator on \, $E_l\, H_{_{ \Pi}} \, .$

Because of the same reason $\, E_l\, \Pi (g)\, (I\, -\, E_l\, )\, =\, 0\, $
 or $\, E_l\, \Pi (g)\, = \,  E_l\, \Pi (g)\,E_l \, .$
\medskip

If we use the similar  reasoning for $g^{^{-1}} \,  ,$ we'll get : \ $
E_l\, \Pi (g^{^{-1}})\, = \,  E_l\, \Pi (g^{^{-1}})\,E_l \, .$
Therefore, \ $  \Pi (g)\,E_l \, = \, E_l \, \Pi (g)\,  .$
The proof of theorem 5.1 is complete.

\vskip 8pt

Let $\xi \,$ be a cyclic vector for the representation
 $ \Pi \  ,\ \Pi (u) \xi =\xi
\, $ for all \, $ u \in U(G(p,\infty ))\, , \, p\, >\, {\bf r}(\Pi )\, $
 and
$q\, \in \, N\, .$

\vskip 8pt
{\bf Proposition 5.2.} \ \ {\it Put \, ${\xi}_q\, = \, \Pi (\sigma )
 \xi \,  ,$
\,  where $\sigma $ equals to $\sigma_q^{(n)}$ or to $ \sigma_q \,\ , \, \,
 H_q \, = \,[\Pi (G(q,\infty)) {\xi}_q\, ]\ . $  ôhen

\hskip 20pt  a)\ \ $  H_q \, = \, \Pi (\sigma )\,  H_{_{ \Pi}}\, \, ; $

\hskip 20pt   b)\ \ representation \ $( \Pi \, ,\, G\, ,\, H_{_{ \Pi}}\,  )
 $ \ is unitary

equivalent to the representation defined by the

chain of mappings:
$$G\, \stackrel{i_q}{\longrightarrow} \, G(q, \infty )\,   \stackrel
 {\Pi }{\longrightarrow} \,
( \Pi (G(q,\infty))\, ,\, H_q \, , \, {\xi}_q\, )\ ,$$

where  \, $ i_q(g) \, $ for $g\in G $ is defined by the correlations

\ ${\sigma}^*\, i_q(g) \,  {\sigma}\, = \, g\, ,\, i_q(g)e_l \, =
 \, e_l \, \;$ for $\; l=\, 1,\, 2,\, 3, \, \ldots  ,\, q \, .$}
\medskip

{\bf Proof. }\  We note, that \ $[\Pi ({\sigma}^*) \Pi (G(q,\, \infty))
\Pi (\sigma )\, \xi ]\, = \, [ \Pi (G) \, \xi ]\, = \,  H_{_{ \Pi}}\, .
 $ \ Therefore,
\, $ \Pi ({\sigma}^*)\, [\, \Pi (G(q,\, \infty))\, {\xi}_q ] \, =
 \, H_{_{ \Pi}}\, $
 or \ $ E_q\, [\Pi (G(q,\, \infty))\, {\xi}_q ]\, =
 \, \Pi ({\sigma}) \,H_{_{ \Pi}}\, \, , $
where \, $ E_q\, = \, \Pi ({\sigma})\,\Pi ({\sigma}^*)\, .$

But by the theorem 5.1 \, $E_q\, \in\, {\Pi (G(q,\, \infty))}^{\prime }\, .$
 Therefore \,
$ \Pi (\sigma) \,  H_{_{ \Pi}}\, = \, E_q\, [\,
 \Pi (G(q,\, \infty))\, {\xi}_q ]\, = \,
[\, \Pi (G(q,\, \infty))\, {\xi}_q ]\, = \, H_q \, .$ Thus correlation
(a)\, is proved.
\medskip

Statement (b)\, follows from this chain of equalities :
\newline
\medskip

$(\Pi (g)\xi\, ,\, \xi )\, = \, (\Pi ({\sigma }^*\, i_q(g)\,
 {\sigma })\, \xi\, , \, \xi )\, =\,
(\Pi ({\sigma }^*) \, \Pi (i_q(g) )\,
 \Pi ({\sigma })\, \xi \, ,\, \xi )\, =\,
(\Pi (i_q(g))\, {\xi}_q \, ,\, {\xi}_q)\, .$
\vskip 9pt

The following statements 5.3-- 5.7 \, refer only to group $GL(\infty)\, .$
Let $q\, = \, 2(p+1)\, ,$ where $p\geq {\bf r} (\Pi)\, .$
\medskip

{\bf Lemma 5.3.} \ \ {\it In $\,  H_{_{ \Pi}}\, $ there exists a set
 of vectors \,
$ \{ {\xi}_i \}\, (i\in N\ ,\ {\xi}_1 \, = \, \xi )\, $with the properties :

a) \ the subspaces\, $H_q(i)\, = \, [\Pi (G_q^I)\, {\xi}_i ]\ (q=2(p+1))\, $

are orthogonal in pairs for the different values of  index $i$

and \ $\oplus_i  H_q(i)\, = \, \,  H_{_{ \Pi}}\, ;$

b)\ \ $\Pi ({\xi}_i)\, =\, { \xi }_i \, $ for all \, $i \in N\, $ and $\,
u \in U(GL(q,\infty ))\, .$}
\medskip

{\bf Proof. } \  Let\ \, $\{ g_k\}  \ \, (k\in N)\ $ be a dense subset in
\, $GL(q) $\, ;\ \, ${\xi}_1\, , \,   {\xi}_2\, ,\, \ldots \, ,\,
 {\xi}_n \, $
are the unit vectors from $H_{_{ \Pi}}\, ,\ $ for which the subspaces
\, $H_q(i)\, = \, [\Pi (G_q^I)\, {\xi}_i ]\ \, (1\leq i \leq n)\ $
 are orthogonal in pairs.
\ If \ $k(n)\, =\, \min\{ k\, :\, \Pi(g_k)\xi\,
 \notin \ \oplus_{i=1}^n\, H_q(i)\,\ \} $ and
\, $ H_{_{ \Pi}}(n)\, =\, \{H_{_{ \Pi}}\, -\,\{\oplus_{i=1}^n\, H_q(i)\,\}\}\
 \neq 0 \, ,$ then put
 ${\xi }_{n+1}\, =\, {\frac{\displaystyle P_n\,
 \Pi (g_{k(n)})\,\xi }{\displaystyle || P_n\, \Pi (g_{k(n)})\,\xi ||}}\, ,$
where $P_n \, $ is an orthoprojection on  $ H_{_{ \Pi}}(n)\, .$

If we continue this process, we will get the system of vectors \,
 $\{{\xi}_i \}\ (i \in N)\, $
 (possibly the finite one) such that\ $\Pi(g_k)\xi\, \in \
 \oplus_{i=1}^{\infty}\, H_q(i)\,\ $
for all  $k\in N$ and \ $\Pi (u)\xi_i\, =\, \xi_i \,\,
 \forall\, u\in U(q,\infty)\, .$

Now this statement follows from lemma 4.3.
\medskip
\vskip 6pt

Let \, $E_q(i)\, $ be the orthoprojection on $H_q(i) \ (H_q(i)\, =
 \, E_q(i)\, H_{_{ \Pi}})\, .$
\, From lemma 5.3 it follows that exists $i\, ,$
 for which $\, E_q(i)\xi_q\,\neq 0\, .$
 If \, $E_q(i)E_q\, =\, V_q(i)\,{[E_qE_q(i)E_q]}^{1\over 2}\, $ is a polar
 decomposition of
$\, E_q(i)E_q\, ,\, e_q(i)\, =\, {V_q^*}(i)\, V_q(i)\,\leq \, E_q\, ,\,
f_q(i)\, =\, V_q(i)\, {V^*}_q(i)\, \leq \, E_q(i)\, ,$ then from theo\-rem 5.1
 and lemma 5.3
we get, that $\,  V_q(i)\,\in {(\Pi(GL(q,\infty )))}^{\prime}\, .$
\medskip

Since $\, {(\Pi (GL(\infty)))}^{\prime\prime}\, =
 \,{( {(\Pi (GL(\infty)))}^{\prime })}^{\prime}\, $
is a factor, then by the proposition 5.2 $w^*\, -$algebra
 ${(\Pi (GL(q,\infty)))}^{\prime\prime}\, $\   of operators,
that act in $H_q\, ,$ is also factor. Because of this, \
accounting theo\-rem 5.1, we get, that the next chain of the mappings
\newline
\medskip

$\, a\in \, {(\Pi (GL(q,\infty)))}^{\prime\prime}\, E_q\, \longrightarrow\,
 e_q(i)\,a\, \in\,  {(\Pi (GL(q,\infty)))}^{\prime\prime}\, e_q(i)\,
 \longrightarrow\,
\newline \longrightarrow \, V_q(i)\, e_q(i)\, a\, {V_q(i)}^* \, =
 \, f_q(i)\,a\,  \in\,  {(\Pi (GL(q,\infty)))}^{\prime\prime}\, f_q(i)\, $
\newline
\medskip

is an isomorphism $\, {(\Pi (GL(q,\infty)))}^{\prime\prime}\, E_q\, =
 \newline =\, E_q\,{(\Pi (GL(q,\infty)))}^{\prime\prime}\, E_q\, $ on
$  f_q(i)\, {(\Pi (GL(q,\infty)))}^{\prime\prime}\, f_q(i)\, =
 \newline  f_q(i)\, {(\Pi (GL(q,\infty)))}^{\prime\prime}\, \subset\,
 {(\Pi (GL(q,\infty)))}^{\prime\prime}\, E_q(i)\, .$
\medskip

Therefore, by the proposition 5.2 the class of unitary equivalence of the
 representation $\Pi $
of group $GL(\infty)\, $ is determined up to the multiplies by the
 restriction of $\Pi $ to
$GL(q,\infty) \, ,$ that acts in some invariant subspace
$\, f_q(i)H_{_{\Pi}}\, \subset \, H_q(i)\, $ with a cyclic vector
 \ $\xi_q(i)\, =\,  V_q(i)\, \xi_q\, .$
\medskip

Let
\begin{equation}
(\Pi (G_q^I),H_q(i),\xi_i )\, =\,\int\limits_S
({\Pi}_s (G_q^I) ,H_q(i,s), \xi_i(s) )\, d\mu (s) \
\end{equation}
be  a decomposition of the restriction of representation $\Pi\, $  to group
 $G_q^I\, ,$
that acts in $H_q(i)\, ,$ into a direct integral of the irreducible spherical
 representations,
that corresponds to center  $C_i $ of algebra
 \ ${(\Pi (G_q^I))}^{\prime\prime}\, ,$ where
$S \, $ is a spectrum of $C_i\, ,\, \mu\, $ is a
 probability measure on $\, S \, .$

 \vskip 6pt
Using the classification of the spherical representations of group
 $GL(\infty)\, ,$
the proposition 3.2, theorem 4.1 and the fact that
 \,$q>2({\bf r }(\Pi)+1)\, ,$
we can prove the following statement.
\medskip

{\bf Lemma 5.4 }\ \ {\it There exist selfadjoint matrix $A$ with a size
 $ \,{\bf r }(\Pi)\times {\bf r }(\Pi)\, $ and  a real number $\beta $
such that
$$({\Pi}_s(g)\xi_i(s)\, ,\, \xi_i(s))\, {||\xi_i(s) ||}^{-1}\, =
 \,$$$$ {\det(|g|)}^{i\beta }\,
\det[I_{{\bf r }(\Pi)}\otimes \cosh (\ln{|g|})\,-\, 2i A \otimes
 \sinh (\ln{|g|})] $$
 $\, \forall\ \  g \in\, GL(q,\infty )\, $ and $\mu\, -$ almost all (a.a.)\,
 $s\in S\, .$
Moreover, for $\mu\, -$ a.a. \,\  $s\in S\,\ \ \newline [\Pi_s (G_q^I)\xi_i(s) ]\, =
 \, [\Pi_s(N_q)\xi_i(s) ]\, ,$
where $\,N_q\,$ consists of matrices of the form \
 $\pmatrix{I_q&0\cr *&I\cr}\ .$}
\medskip

From this and from the complete classification of the spherical
 representations of groups of motions
$G_n^i\, $ (see\, \cite{N2} ) it follows
\medskip

{\bf Lemma 5.5.}\ \ {\it  There're exist : isometry $V\, ,$ that maps
 $\, H_q(i) \, $
to $\, L^2(S,\mu )\otimes L^2({\Lambda}_{_{{\bf r }(\Pi)}}\, ,\,
 {\nu}_{_{{\bf r }(\Pi)}})\, ,$
and $\mu \, -$ measurable mapping  $z\,$ from $\, S\, $ to the set of all
 the complex $\ q\times {\bf r }(\Pi)\,$ matrices
such that the action of the operators
$\ {\tilde\Pi}_{_{\scriptstyle A}}(g)\, =
 \, V\Pi (g) E_q(i) V^{-1}\ (g\in G_q^I) \ $
is determined by correlations :

\begin{equation}
({\tilde\Pi}_{_{\scriptstyle A}}(g)\eta )(s)\, =\,
 {\Pi}_{_{\scriptstyle{Az(s)}}}\eta (s)\ \
\end{equation}\ \ (see (6))\, ,
where $\,  \eta (s)\, \in\,
 L^2({\Lambda}_{_{{\bf r }(\Pi)}}\, ,\, {\nu}_{_{{\bf r }(\Pi)}})\, $
for all $\, s\in S\ .$

Moreover, there exists a \,$\mu\, -$ measurable mapping
 $\, f_q(i,\, .)\, $  from $S\, $ to the set of
the orthoprojections of $\, w^*\, -$algebra
 $\, {(\Pi_{_{\scriptstyle{Az(s)}}}
 (GL(q,\infty)))}^{\prime}\, \subset\, B( L^2({\Lambda}_{_{{\bf r }(\Pi)}}\, ,
 \, {\nu}_{_{{\bf r }(\Pi)}})\,)\, ,$
for which\newline  $\ (Vf_q(i)V^{-1}\eta )(s )\, =\, f_q(i,s)\eta (s)\ .$}
\vskip 8pt
\medskip

For $\, u \in U({\bf r }(\Pi ) , A)\, =\, \{u \in U({\bf r }(\Pi )) :
 \, [u\, ,\, A]\, =\, 0\}\,$
we define operator $\, \tau (u)\, $ in
 $\  L^2({\Lambda}_{_{{\bf r }(\Pi)}}\, ,\, {\nu}_{_{{\bf r }(\Pi)}})\, $
by
\begin{equation}
(\tau (u)\xi )(\lambda)\, =\, \xi (u^* \lambda)\ .
\end{equation}

Denote by $\ \tilde{\tau }(u)\, $ a natural extension of $\, \tau\, $ to
$\ L^2(S,\mu )\otimes
 L^2({\Lambda}_{_{{\bf r }(\Pi)}}\, ,\, {\nu}_{_{{\bf r }(\Pi)}})\, .$
 Namely, \newline $\ \tilde{\tau }(u)\, =\, I\otimes {\tau }(u)\, .$
\medskip

If $\kappa \, $ is irreducible representation of group $\, U({\bf r }(\Pi ) , A)\, ,\,
{\kappa}_{ik}\, $
is a matrix element of operator $\kappa (g)\, ,$ then

\medskip

$\ P_i^{\kappa }\, =\, \dim ({\kappa })\int\limits_{U({\bf r }(\Pi ) , A)}
{\kappa}_{ii}(v)\tau (v)\, d v\ $
 is an orthoprojection from $w^* -$algebra\newline
 $\ {({\Pi}_{_{\scriptstyle{Az}}}
 (GL(q,\infty)))}^{\prime} \, $
(see (6)), that acts in $\
 L^2({\Lambda}_{_{{\bf r }(\Pi)}}\, ,\, {\nu}_{_{{\bf r }(\Pi)}})\, .$
\medskip

Obviously  operators   $\, {\tilde P}_i^{\kappa }\, =\, I\otimes
 P_i^{\kappa }\, $
and $ \ \tilde{\tau }(u)\, \, \, \forall\  \ u\in U({\bf r }(\Pi ) , A)\, $
belong to \newline $\ {{\tilde\Pi}_{_{\scriptstyle A}}(GL(q,\infty))}^{\prime} \ .$
\medskip

Without loss of generality suppose, that for ${\kappa }$ and $i$ introduced
 above
$\,  {\tilde P}_i^{\kappa }\, Vf_q(i)V^{-1}\, \neq \, 0  \, .$

\vskip 8pt

Let $\ w|Vf_q^i V^{-1} {\tilde P}_i^{\kappa}Vf_q^i V^{-1}|^{1\over2}\, $
 be a  polar decomposition
of operator $\  {\tilde P}_i^{\kappa}Vf_q(i) V^{-1} \, .$
 Since $\ {({\tilde\Pi}_{_{\scriptstyle A}}
 (GL(q,\infty)))}^{\prime\prime}Vf_q(i)V^{-1}\ $ is a factor ,
 orthoprojection $\, w^*w\leq Vf_q(i) V^{-1} \, ,\, w\, ,\,
\newline Vf_q(i) V^{-1}\, \in\,
{({\tilde\Pi}_{_{\scriptstyle A}}(GL(q,\infty)))}^{\prime}\ ,$ then, using
the statement of theorem 0.1  and the fact that\newline orthoprojection
 $\,  P_i^{\kappa}\, $is a minimal in $\
 {(\Pi_{_{\scriptstyle{Az(s)}}} (GL(q,\infty)))}^{\prime}\, \subset\,
 B( L^2({\Lambda}_{_{{\bf r }(\Pi)}}\, ,\, {\nu}_{_{{\bf r }(\Pi)}})\,)\, \ $,

 we get:

{\it i)}\ the restriction  of $\, {\tilde\Pi}_{_{\scriptstyle A}}\, $ to
 $\, GL(q,\infty)\, ,$ that acts in

$\, Vf_q(i)V^{-1}(\, L^2(S,\mu )\otimes L^2({\Lambda}_{_{{\bf r }(\Pi)}}\, ,
 \, {\nu}_{_{{\bf r }(\Pi)}})\, )\, ,$  is multiple

by the representation of subgroup  $\, GL(q,\infty)\, $ in

$\, ww^*(\, L^2(S,\mu )\otimes
 L^2({\Lambda}_{_{{\bf r }(\Pi)}}\, ,\, {\nu}_{_{{\bf r }(\Pi)}})\, )\,)\, $
  defined

by correlation $\, w {\tilde\Pi}_{_{\scriptstyle A}}(g)w^*\, =
 \, {\tilde\Pi}_{_{\scriptstyle A}}(g)ww^*\, ;$

{\it ii) }\  if $\, f(s) \,\,\mbox{ is}\ \, \mu\, -$measurable field of
 orthoprojections,

that determines orthoprojection $\,  ww^*\, ,$ then there exists

 subset $\, S^{\prime} \subset S\, $ of the positive measure such that
\medskip

 $ f(s) = \left\{
\begin{array}{rl}
P_i^{\kappa}, & \mbox{if}\ \  s \in S^{\prime} \\
0 ,           & \mbox{if}\ \  s \not\in S^{\prime}
\end{array} \right. $

\vskip 8pt
The facts above are summed by the following statement.
\medskip

{\bf Proposition 5.6.}\ \ {\it In the commutant of representation
 $\, (\Pi\, ,\,GL(q,\infty)\, ,\, H_q\, ) \, $
(see proposition 5.2) there exists an orthoprojection $f\, $ such, that
$\, (\Pi\, ,\,GL(q,\infty)\, ,\, fH_q\, ) $  is irreducible and unitary
equivalent to the restriction of $\, {\Pi}_{_{\scriptstyle{Az}}} \, $
 to group $\, GL(q,\infty))\, ,$
that acts in $\, P_i^{\kappa}
 L^2({\Lambda}_{_{{\bf r }(\Pi)}}\, ,\, {\nu}_{_{{\bf r }(\Pi)}})\, . \, $ }

\vskip 8pt
From this and from proposition 5.2 (b) it follows the main classification

{\bf Theorem 5.7.} \ { \it An arbitrary admissible factor--representation
 $\, \Pi \,$
of group $\, GL(\infty)\, $ has a type $\, I\, $ and it is multiple
by representation  $\,{\Pi}_{_{\scriptstyle{Az}}} \,$
  of group $G_n^I \,$ for $n=0\, ,\, z=0 \ (G_o^I\, =\, GL(\infty)\, )\, ,$
 that acts in
$\,  P_i^{\kappa}
 L^2({\Lambda}_{_{{\bf r }(\Pi)}}\, ,\, {\nu}_{_{{\bf r }(\Pi)}})\, \,  ,$
for some selfadjoint matrix $\, A\, $
  of size $\, {\bf r }(\Pi)\times {\bf r }(\Pi)\, $  and orthoprojection  $\,
 P_i^{\kappa} \,  ,$
where $\kappa \, $ is irreducible representation of $\, U({\bf r }(\Pi),A)\, ,$
(see (6)) . }
\bigskip

\bigskip

\vskip 8pt
A statements analogous to statements (5.3) -- (5.6), from which follows
 theorem 5.7,
might be modified for group $G_n^I \ .$
\medskip

Let $\,\Pi\, $ be an admissible factor--representation of group
 $G_n^I\,\,  (n\geq 1)\, , \, \xi\, \;$ is
a cyclic vector for $\Pi \, ,\, {\bf r }(\Pi)\, \,$ is a rang of
 $\Pi\, $ (see the definition 3.5),
$\, \Pi (u)\xi\, = \, \xi\,$ for all $\, u\in \, U(G(p+n,\infty)) \, ,\,
q=n+2(p+1)\, ,\, p\geq  {\bf r }(\Pi)\, .$
\medskip

The following statement is analogous to lemma 5.3 for group $\, G_n^I\ .$

{\bf Lemma 5.8. } {\it \  In $\, H_{_{\Pi}}\, $ there exists a set
of unit vectors $\, \xi_i \,\, (i\in N)\ (\xi_1\, =\, 1)  \   $
 with the properties :

a) subspaces $\, H_q(i)\, =\, [\Pi (G_q^I)\xi_i ] \ (q=n+2(p+1))  \,  $
 are orthogonal

in pairs   for different $i$  and $\, \oplus_{i\in N} H_q(i)\, =\,
 H_{_{\Pi}}\, ;$

b)\  $\Pi (u)\xi_i\, =\, \xi_i\, $ for all $i\ $ and
 $\, u\in U(G(q,\infty))\, .$ }
\medskip

{\bf Proof} might be done by using the statement of lemma 4.4
 and analogous to
the ground of lemma 5.3.

\vskip 8pt

If $\, E_q\, ,\, \xi_q \, $ are the same as in proposition 5.2, then
repeating the reasoning we've done for $\, GL(\infty)\, ,$
 let's find the orthoprojection
$\, f_q(i)\,\in {(\Pi (G(q,\infty)) )}^{\prime }  $ such that
 $\, f_q(i) H_{_{\Pi}}\, \subset \, H_q(i)\, $ and
$ (\Pi\, ,\,  G(q,\infty)\, ,\, H_q )\ $ is multiple by
 $\ (\Pi\, ,\,  G(q,\infty)\, ,\, f_q(i)H_q(i)\, .$

\vskip 6pt
Hence accounting  theorem  2.1, statements 3.2, 4.1 and the decomposition
\newline
$(\Pi (G_q^I)\, ,\,H_q(i)\, ,\,\xi_i)\, =\,
\int\limits_S (\Pi_s (G_q^I)\, ,\, H_q(i,s)\, ,\, \xi_i (s)\, )\, d\, \mu(s)
 \ ,$
where the objects we've met are described in  (11), we get the statement
 analogous to lemma 5.5 for $\, G_n^I\, .$
\medskip

{\bf Lemma 5.9.} {\it  There exist : s.a.  $\, {\bf r }(\Pi)\times
 {\bf r }(\Pi)\, -$
matrix $\, A\, ,\, \mu -$ measurable mapping $z\, $ from $\, S\,$ to the set
 of matrices of
$q=n+2(p+1)$ rows and ${\bf r }(\Pi)$ columns of a form $\, z(s)\, =\,
 \pmatrix{z_n\cr *\cr }\, ,
  $ where $\, z_n\, $ is independent of $s\, $ matrix  of height $n \, ;
 $  isometry
$\, V\, , $ that maps $\, H_q(i)\, $ to $\, L^2(S,\mu)\otimes
 L^2({\Lambda}_{_{{\bf r }(\Pi)}}\, ,\, {\nu}_{_{{\bf r }(\Pi)}})\, $
such that  operators $\tilde{\Pi}_{_{\scriptstyle {Az}}}(g)\, =
 \, V{\Pi}_i(g)V^{-1}\,  \, (g\in G_q^I) \, ,$  where ${\Pi}_i\, $ is
a restriction of  $\Pi\, $ to $\, H_q(i)\, ,$
are defined in $ \, \, L^2(S,\mu)\otimes
 L^2({\Lambda}_{_{{\bf r }(\Pi)}}\, ,\, {\nu}_{_{{\bf r }(\Pi)}})\, $
by the correlations (6)\, and\, (12).}

\vskip 6pt

\vskip 6pt
Let $\,U( {\bf r }(\Pi),A,z_n) =\{ u\in U( {\bf r }(\Pi),A)\, :\,
 z_n\, u\, =\, z_n \}\,  .$
As  before we define by  the irreducible representation $ \kappa\, $
 of group $\,U( {\bf r }(\Pi),A,z_n)\, $
the orthoprojection $\, P_i^{\kappa}\, .$

\vskip 6pt
The reasoning for group $GL(\infty )\, ,$ introducing in this chapter with a
  slight
change might be carried over to the case of group $\, G_n^I\, $ and
 they're leading to the following statement .
\medskip

{\bf Theorem 5.10.} {\it  For the admissible factor--representation
 $\,\Pi \,$ of group   $\, G_n^I\, $
there exist : s.a.  $\, {\bf r }(\Pi)\times {\bf r }(\Pi)\, -$
matrix $\, A\, ,\, n\times {\bf r }(\Pi)\, -$ matrix $z\, $ and
 orthoprojection
$\, P_i^{\kappa}\, $ such that $\Pi $ is multiple by  restriction of
 representation
${\Pi}_{_{\scriptstyle {Az}}}\, $ (see (6)) \ to subspace
$\, P_i^{\kappa}\,
 L^2({\Lambda}_{_{{\bf r }(\Pi)}}\, ,\, {\nu}_{_{{\bf r }(\Pi)}})\, .$}

\vskip 27pt
\begin{center}
{\Large
\S 6  .   A Description Of The Admissible Representations   \newline
        \hskip 22pt     Of Groups $\ Sp(2\infty)\ $ \ And \ $O(2\infty) $}

\end{center}

\vskip 8pt

The main result of this chapter is theorem 6.15, that establishes a
completeness of the collection
of the admissible factor--representations, that was built in \S 1, for
groups $Sp(2\infty)\,$ and $O(2\infty) .$
    We suppose that  $G$ coincides with the one of the groups
$Sp(2\infty) $ or $O(2\infty) .$ Otherwise we'll specialty notice it.

Let ${\Theta}_1^{(n)}\, $ be a subgroup of $\, {\Theta}^{(n)} \, ,$ that
consists of
matrices of the form $\,  {\theta}^{(n)}(x,0)  \, ,\ GN_n\, =\,
G_d(n,\infty )\, {\Theta}_1^{(n)}\,  .$
Notice, that subgroup $\ GN_n\, $ is naturally isomorphic to the group of
motions $\, G_n^I \, .$
\medskip

{\bf Proposition 6.1.} {\it If  $\,\Pi \, $ is an admissible factor--
representation of
 group $G\, ,$ of the rang $\, {\bf r }(\Pi)\, $ that acts in a Hilbert
space
$\, H_{_{\Pi}} \, ,\, \eta $ is a nonzero vector from  $\, H_{_{\Pi}} \, $
and $n \geq  {\bf r }(\Pi)\, ,$
then the following properties are true :
\vskip 5pt

{\sl i)}\  in a subspace $\, H_{\eta }\, =\, [\Pi (GN_n)\eta ]\ $
there exists

a nonzero vector ${\eta }_{_{U}}\, ,$ for which  $\Pi (u){\eta }_{_{U}}\, =
\, {\eta }_{_{U}}  \, $

$\forall \ \, u\in U( GN_n)\ \, ;     $
\vskip 6pt

{\sl ii)} if $\ (\Pi (GN_n), H_{\eta },\eta)=\int\limits_S(\Pi_s (GN_n),
H_{\eta }(s),\eta(s))  d\mu(s) \, -$

is a decomposition of  $\  (\Pi ,GN_n,H_{\eta })\ $ into a direct integral

of the factor--representations corresponding to the center of

${( \Pi (GN_n))}^{\prime\prime}\, ,$ then for $\mu -$ a.a.
$s\in S\ \ (\Pi_s , GN_n, H_{\eta }(s)) \, $

is multiple by $({\Pi}_{_{\scriptstyle {Az(s)}}}\, ,\, GN_n\, ,\,
P_i^{\kappa}L^2({\Lambda}_{_{{\bf r }(\Pi)}}\, ,\, {\nu}_{_{{\bf r }(\Pi)}}))
\, $

(see theorem 5.10 and (6)). The rang of  matrix $z(s)\, $ being

 equals to ${\bf r }(\Pi)\, $ for $\mu\, -$a.a. $\, s\in S .$}
\medskip

{\bf Proof } is based on the theorem 4.1, the classification statements
5.7, 5.10 and the analyses of the concrete realization of  representations
of group
$\, G_d\, ,$ that is isomorphic to $\, GL(\infty)\, .$
\vskip 8pt

{\bf Remark 6.2 }\ \ {\sl Let $\Pi\, $ be the same as in a proposition 6.1,
$\xi\, $ is a cyclic vector for $\Pi\, \, (\, H_{_{\Pi}}\, =
\, [\Pi (G)\xi ] )\, .$
ôhen for $\Pi $ and $G\, $  the theorem 5.1 and the proposition 5.2 are
true. Namely,
a class of the unitary equivalence of representation $\, \Pi\,$
of group $G\, $
is uniquely determined by restriction of $\Pi \,$  to $G(q,\infty)\, ,$
that acts in
$\, H_q\, =\, [\Pi (q,\infty){\xi }_q]\, ,$ where $\, \xi_q\, =
\, \Pi (\sigma_q)\xi\,
\ \  (\sigma_qe_i\, =\, e_{i+(sign i)q})  $ \
(see proposition 5.2).}
\vskip 6pt

We denote by $E_q\, $ the projection of $H_{_{\Pi}} \, $ to $H_q\, .$
\medskip

From the proposition 6.1 there follows an important

{\bf Lemma 6.3.}\ \ {\it Let $\Pi\, $ be the same as in conditions of
statements
6.1 -- 6.2, $q\, =\, {\bf r }(\Pi) .$ There exists $\Pi (U(GK_q))\, -$
fixed  unit vector
 ${\eta}_q^{^{(U)}}\, \in \, H_{_{\Pi}}\, $ with a property:
if $E_q^{^{(U)}}\, $
is an orthoprojection of $\, H_{_{\Pi}}\, $ to
$\, [ \Pi (GK_q){\eta}_q^{^{(U)}}\,]\, ,$
then $E_q^{^{(U)}}\xi_q\, \neq 0 \, .$}
\medskip

{\bf Proof. }\ \ As  $ {\eta}_q^{^{(U)}}\,$ we should take
the unit  $\Pi (U(GN_q))\, -$ fixed vector from  subspace
$H_{\xi_q}\, =[\Pi (GN_q)\xi_q ]\ ,$
that exists by proposition 6.1 ({\it i}). It satisfied to all the conditions
of
the lemma, but it may be not  $\Pi (U(GK_q))\, -$ fixed.
The last property follows from
the complete description of the structure of  an admissible
representations of group $U(GK_q)\, ,$
that was found by á.á. Kirillov in \cite{K},\,  and  from the fact that
vector  $ {\eta}_q^{^{(U)}}\,  $ is   $\Pi (U(GN_q))\, -$ fixed.

Let $\, v_q\,{|E_q E_q^{^{(U)}} E_q|}^{1\over 2 }\, =
\, E_q^{^{(U)}} E_q \, $
be a polar decomposition of operator $\, E_q^{^{(U)}} E_q \, , \,
{ v_q}^*v_q\, =\, e_q\leq E_q\, ,\
v_q { v_q}^*\, =\,   e_q^{^{(U)}}\leq  E_q^{^{(U)}}\ .$
\medskip

 Diminishing, if we need, orthoprojection $\,e_q \, $ and accounting the
fact that
$ v_q \in {\Pi (G(q,\infty )) }^{\prime}\, ,$ Ánd
$ {\Pi (G(q,\infty )) }^{\prime\prime}\, $
is a factor in $H_q\, ,$ and using the  statements 6.1 -- 6.2 we get

{\bf Proposition 6.4.} {\it  Let $\Pi $ be the same as in
statements 6.1 -- 6.3.
The representation $\, (\Pi \, ,\,  G(q,\infty )\, ,\, H_q)\, ,$ that is a
restriction
of $\Pi\, $ to the group $\, G(q,\infty )\, $ and that acts in
$\, H_q\, =\, [\Pi (G(q,\infty ))\xi_q]\, $
is multiple by  $\, (\Pi \, ,\,  G(q,\infty )\, ,\,
e_q^{^{(U)}}H_{_{\Pi}})\, ,$
 where $\,  e_q^{^{(U)}}H_{_{\Pi}})\, \subset \, [\Pi (GK_q)v_q\xi_q ]\,
\subset \, [\Pi (GK_q){\eta}_q^{^{(U)}}]\, =
\, H_q^{^{(U)}}\, ,\, {\eta}_q^{^{(U)}}\, $
is a unit $\, \Pi (U(GK_q))\, -$ fixed vector.}
 \medskip

{\bf Proof } follows from the correlation $\, v_q \Pi (g)v_q^*\, =
\, v_q e_q\Pi (g)v_q^*\, =\, v_q e_qv_q^*\Pi (g)\, =
\, e_q^{^{(U)}}\Pi (g) \,\newline \forall g\in G(q,\infty )\, .$

Let's begin the analysis of restriction of $\Pi \, $ to $GK_q\, ,$
that acts in $\, H_q^{^{(U)}}\, .$
\medskip

{\bf Proposition 6.5.}\ {\it  Let  $C(M)\, $ be  a center of $w^*\, -$
algebra $M\, ,\, \Pi\, ,\, {\eta}_q^{^{(U)}}\, $
are the same as in the condition of proposition 6.4, $q\, =
\, {\bf r}(\Pi)\, .$
 ôhen $ C({\Pi (GN_q)}^{\prime\prime})\, \subset\,
C({\Pi (GK_q)}^{\prime\prime})\, .$}
\medskip

{\bf Proof}\ \ Let $ {\Theta }^{(n,q)}\, =\, {\Theta }^{(n)}\cap
{\Theta }^{(q)}
\ \, (n\geq q)\, ,\   {\Theta }_1^{(n,q)}\, =\,
{\Theta }_1^{(n)}\cap  {\Theta }^{(n,q)}\, ,\,
\ GK_{n,q}\ (GN_{n,q})\ $ is generated by $\ G(n,\infty)\, ,\
\Delta_q\, $\, and\, $ {\Theta }^{(n,q)}\, \ (G_d(n,\infty )\, $ and
$\, {\Theta }_1^{(n,q)}\, )\ .$
Obviously $  GK_{q,q}\, =\,  GK_{q}\, ,\   GN_{q,q}\, =\,  GK_{q}\ .$
\medskip

If we prove  the correlations
\begin{equation}
 C({\Pi (GK_q)}^{\prime\prime})E_q^{^{(U)}}\, =\,\bigcap\limits_{n\geq q}
C({\Pi (GK_{n.q})}^{\prime\prime})E_q^{^{(U)}}\, ,
\end{equation}
\begin{equation}
 C({\Pi (GN_q)}^{\prime\prime})E_q^{^{(U)}}\, =\,\bigcap\limits_{n\geq q}
C({\Pi (GN_{n.q})}^{\prime\prime})E_q^{^{(U)}}\, ,
\end{equation}
then from the fact that $GN_q\, \subset GK_q\, ,$ our statement will
follow.
\medskip

Let $\, c\in \, C({\Pi (GK_q)}^{\prime\prime })E_q^{^{(U)}}\, .$
There exist the collections
of the elements
$\, {\bigg{\{} g_{_{\displaystyle i k_i}}\bigg{\}} }_{_{\displaystyle i k_i }
}^{N_{_{\scriptstyle i}}}\, \subset\, \newline \,  GK_{q}\cap G(i)\, $
and the complex numbers   $\, \, {\bigg{\{}
c_{_{\displaystyle ik_i}}\bigg{\}} }_{_{\displaystyle ik_i}}^{
N_{_{\scriptstyle i}}}\,  \, (i\in N) $
 with the properties:
$$\bigg|\bigg| \sum_{{\displaystyle k_i}=1}^{N_{_{\scriptstyle i}}}
c_{_{\displaystyle ik_i}}  \Pi (g_{_{\displaystyle ik_i}})\bigg|
\bigg| <||c||, $$
\begin{equation}
\lim_{i\to\infty}\bigg|\bigg| \sum_{{\displaystyle k_i}=
1}^{N_{_{\scriptstyle i}}} c_{_{\displaystyle ik_i}}
\Pi (g_{_{\displaystyle ik_i}})f-cf\bigg|\bigg|=0 \ \ \forall\ \
f\,\in\,  H_q^{^{(U)}}\, .
\end{equation}

Next we notice, that in $U(G_d(q,\infty)) \, $
there exist a sequence of elements
 $u_i\, \, (i\in N)\ $ with the property
\begin{equation}
u_i\, g_{_{\displaystyle i k_i}} \, u_i^*\ =
\ { g^{\prime }}_{_{\displaystyle i k_i}}\, \in\,
GK_{i,q}\ \ \forall\ i\,\geq\, q\ .
\end{equation}

Hence accounting (16) and the fact that  vector $\, {\eta}_q^{^{(U)}}\, $
is $\Pi (U(G_d(q,\infty)))\, -$fixed
 we get
$$\lim_{i\to\infty}\bigg|\bigg| \sum_{{\displaystyle k_i}=
1}^{N_{_{\scriptstyle i}}} c_{_{\displaystyle ik_i}}
\Pi ({ g^{\prime }}_{_{\displaystyle ik_i}}){\eta}_q^{^{(U)}}-c\,
{\eta}_q^{^{(U)}}\bigg|\bigg|\, =0 \ . $$

This correlation and (17) lead to the following chain of equalities:
$$0\, =\, \lim_{i\to\infty}\bigg|\bigg|\Pi (g)\bigg( \sum_{{\displaystyle
k_i}=1}^{N_{_{\scriptstyle i}}} c_{_{\displaystyle ik_i}}
\Pi ({ g^{\prime }}_{_{\displaystyle ik_i}}){\eta}_q^{^{(U)}}-
c\, {\eta}_q^{^{(U)}}\bigg)\bigg|\bigg|\, = $$
$$=\,\lim_{i\to\infty}\bigg|\bigg| \sum_{{\displaystyle k_i}=
1}^{N_{_{\scriptstyle i}}} c_{_{\displaystyle ik_i}}
\Pi ({ g^{\prime }}_{_{\displaystyle ik_i}})\Pi (g){\eta}_q^{^{(U)}}-
c\,\Pi (g) {\eta}_q^{^{(U)}}\bigg|\bigg|\, =0 \ , $$
that are true for all $\, g\in GK_q\, .$
\medskip

Hence accounting (16), we get the property (14).
Correlation (15) are proved by the same way .
The proposition 6.5 is proved.
\medskip

{\bf Proposition 6.6.} \ {\it  Let $\, \Pi\, ,\, q\, ,\, {\eta}_q^{^{(U)}}\, ,
\, H_q^{^{(U)}}\, $
are the same as in proposition 6.5,
 $\newline (\Pi, GK_q\, ,\, H_q^{^{(U)}})\, =\, \int\limits_S
(\Pi_s\,,\, GK_q\, ,\, H_q^{^{(U)}}(s) \, ) d\mu(s)\, $
is a decomposition
of $ (\Pi, GK_q\, ,\, H_q^{^{(U)}})\, $ into a direct integral
of the factor--representations,
that corresponds to the center of ${\Pi (GK_q)}^{\prime\prime}
E_q^{^{(U)}}\, ,\, \newline {\eta}_q^{^{(U)}}\, =
\,  \int\limits_S {\eta}_q^{^{(U)}}(s) d\mu(s)\, .$
ôhen

{\sl i)}\ for $\mu\, -$ a. a.  $\, s\in S \, $  subspace
$\{ \eta (s)  \in H_q^{^{(U)}}(s)\, :$

$ \, \Pi_s(u)\eta (s)=\eta (s)\ \forall\
u\in U(GK_q)\,\}\, $ is one-dimensional;

{\sl ii)}\ ${\Pi (GN_q)}^{\prime\prime}\, $ is a factor in
$\, H_q^{^{(U)}}(s)\, ;$

{\sl iii)}\ $ H_q^{^{(U)}}(s)\, =\, [{\Pi_s (GK_q)}{\eta}_q^{^{(U)}}(s)] =
[{\Pi_s (GN_q)}{\eta}_q^{^{(U)}}(s)]\,]\, \, $

for $\mu\, -$a.a. $s\in S .$}
\medskip

{\bf Proof. } \ Since $GK_q\, $is an  Á.Á. group (see definition 3.1),
and $(\Pi_s\,,\, GK_q\, ,\, H_q^{^{(U)}}(s) \, )\, $ is a cyclic factor--
representation,
then the property {\sl (i)} is a corollary of general theorem 3.4.
\medskip

Statement {\sl (ii)} follows from the proposition 6.5.
\medskip

To prove {\sl (iii)} we notice, that  subspace
$H^{\perp }(s) = H_q^{^{(U)}}(s)\ominus
[{\Pi_s (GN_q)}{\eta}_q^{^{(U)}}(s)]\,\,\,\,
\mbox{is}\newline  \Pi_s (GN_q)\, -$invariant.
If $H^{\perp }(s)\,\neq 0\, ,\ P^{\perp }(s)\, $ is an orthoprojection from
$\, H_q^{^{(U)}}(s)\, $ to $\, H^{\perp }(s)\, ,$ then
by the property {\sl (ii)} in  ${\Pi_s (GN_q)}^{\prime}\, $ there exists
a partial isometry $w_s\, ,$ for which $\, w_sw_s^*\, \leq \
P^{\perp }(s)\, $ and $\, w_s{\eta}_q^{^{(U)}}(s)\neq\, 0\, .$
By the construction the vector $\, w_s{\eta}_q^{^{(U)}}(s)\, \mbox{is} \ \
\Pi_s (U(GN_q))\, - $ fixed.
Since $\Pi_s\, $ is an admissible representation for $\mu\, -$a.a.
$s\in S \, ,$ then from the structure
of the admissible representations of groups $\, U(Sp(2\infty))\ $ \ and \
$ U(O(2\infty))   \, ,$ that was completely described
by á.á.Kirillov in  \cite{K}, we get that vector
$\, w_s{\eta}_q^{^{(U)}}(s)\, $   is $\, \Pi_s (U(GK_q))\, - $fixed.
But the last fact contradicts the property {\sl (i)}. Proposition 6.6 is
proved.

The next statement follows from statements 3.2, 6.1, 6.6 and
classification
theorem 5.10.
\medskip

{\bf Proposition 6.7. }\ \, {\it Let  $\Pi\, $ be the same as in
proposition 6.6.
ôhen there exists $\mu\, -$measurable field of the isometries
$\, V_s\ (s\in S)\, ,$
that maps $\, H_q^{^{(U)}}(s)\,$ to $\, L_2(\Lambda_q,\nu_q)\, $
such that
 operators $V(s)\Pi_s(g)V_s^{-1} \, =\, {\widehat\Pi}_s(g)\, $ act
in $\, L_2(\Lambda_q,\nu_q)\, $
for $g\in GN_q \, $ by:
$$ ({\widehat\Pi}_s(g_q)\eta )(\lambda)={\hat{\alpha}_{_A}}(\lambda,g)
\eta (\lambda g)\ , $$
 \begin{equation}
({\widehat\Pi}_s({\theta}^{(q)}(x,0) )\eta )(\lambda)\, =\,
\exp \lbrack i  \Re  Tr (z(s) \lambda h )
\rbrack \eta (\lambda )\, ,
\end{equation}
where $\hat{\alpha}_{_A}\, $ is defined in the condition of
theorem 2.1, $z(s)\,\  \mbox{is} \ \  q\times q\, -$
matrix, $A\, $is s.a. $\,\newline q\times q\, -$matrix and
$q={\bf r}(\Pi ) \, .$

Moreover $A\,$ does not depend on $s\, ,$ the mapping $s\rightarrow z(s) \, \
\mbox{is} \ \ \mu\, -$measurable and
the rang of $z(s)\, $ equals to $q\,$ for $\mu\, -$ a.a. $s\in S\, .$}

\bigskip

Now we should  get ( side by side with (18)) explicitly the form of
the actions of  operators $\, {\widehat\Pi}_s({\theta}^{(q)}(x,y) )\, ,
\, {\widehat\Pi}_s({\gamma }_{_{\scriptstyle u}}^{(q)}(b) )\, $ and
$\ {\widehat\Pi}_s({\gamma }_{_{\scriptstyle o}}^{(q)}(b))\, $ (see \S \,3 ).
\bigskip

The next statement follows from the commutation relations (8) and from
 propositions  6.6 ({\sl iii})-- 6.7.
\medskip

{\bf Proposition 6.8. }\ \ {\it Let $\Pi \, $ satisfies to all the requests
of proposition 6.1, and $\, {\widehat\Pi}_s \, \,\ (s\in S)\, $ are built as
 in proposition 6.7.   There exists $\mu\, -$measurable mapping $h \, $
from $S\,$ to the set of
symmetric (antisymmetric ) for $\, G=Sp(2\infty)\,\, \,
(G=O(2\infty)\,) \,$
$q\times q \, -$matrices such that the operators
$\, {\widehat\Pi}_s(g) \, \ (g\in GK_q)\,$
act in $\, L^2(\Lambda_q ,\nu_q )\, $ by :
\begin{equation}
  ( {\widehat\Pi}_s( {\delta}^{(q)}(a))\eta )(\lambda )=
\exp\{i\,\Re\, Tr(a {h}(s))\}\,\eta (\lambda)\, ,\,
\end{equation}
\begin{equation}
({\widehat\Pi}_s(g_q)\eta )(\lambda)={\hat{\alpha}_{_A}}(\lambda,g)
\eta (\lambda g)\ ,
\end{equation}
\begin{equation}
({\widehat\Pi}_s({\theta}^{(q)}(x,0) )\eta )(\lambda)\, =\,
\exp \lbrack i  \Re  Tr (z(s) \lambda h )
\rbrack \eta (\lambda )\, ,
\end{equation}

$$({\widehat\Pi}_s({\theta}^{(q)}(0,y) )\eta )(\lambda) =
U(y,\lambda , s)\,\cdot $$\begin{equation}
\cdot {\bigg{[}{{d\nu_q(\lambda +
2z^{-1}(s)(y{h}(s))^t)}\over{d\nu_q(\lambda) }}\bigg{]}}^{1\over{2} }
\eta )(\lambda + 2z^{-1}(s)(y{h}(s))^t)\ ,
\end{equation}
\begin{equation}
( {\widehat\Pi}_s({\gamma }_{_{\scriptstyle u}}^{(q)}(b) )\eta )
(\lambda)=
M(b,\lambda , s )\eta (\lambda )\, .
\end{equation}

Moreover $U\, $ and $M\, $ are the unitary scalar functions.}

\vskip 7pt

In the following statement we'll find an explicit form of cocycle $U\,$ and
function $M\, ,$
that are introduced in proposition 6.8.

\vskip 4pt

{\bf Lemma  6.9.}\ \ {\it If $\, {h}(s)\, ,\, z(s)\ \, (s\in S)\, $
are the same as in conditions of \newline propositions 6.7 -- 6.8,
then the following correlations are true:

{\sl i)}$\ U(y,\lambda , s)=\exp\{ 2i \, Tr[{\lambda}^*Az^{-1}(s)
{(yh(s))}^t+$

${({\lambda}^*Az^{-1}(s){(yh(s))}^t)}^*+2{(z^{-1}(s){(yh(s))}^t)}^*
Az^{-1}(s){(yh(s))}^t]\}\ $;

{\sl ii)} \ $\ {h}(s)\ $ for $\mu \, -$a.a. $s\in S\ $ is an  invertible;

 {\sl iii)}\ $\ M(b,\lambda ,s)=\exp\{-i \, \Re\,
Tr[R(s)\lambda b{\lambda}^t]\}\ ,$

where $\, R(s)={1\over4}z^t(s)h^{-1}(s)z(s) \, .$}
\medskip

{\bf Proof.}\ \ In the  following reasoning we omit some standard details
of metric character.

Let $Y(\infty , q)\, $ be a set of all local nonzero  $(\infty\times q\, )-$
matrices
of infinite number of rows and $q\, $ columns. If $y\in Y(\infty , q)\, ,$
then we denote by $y(k)\,$ a matrix from $Y(\infty , q)\, , $ that has only
one
 nonzero $k-$th  row, which coincides with the $k\, -$th row of matrix $y\, .$
Since
operator $ {\widehat\Pi}_s({\theta}^{(q)}(0,y(k)) )\, $  commutes with  all
the
${\widehat\Pi}_s(u) \, , $ where $u\in U(G_d(q,\infty ))\, $ and satisfies
to the
correlation $ue_m=e_m\, $ for $m=\pm k\, ,$ then it's easy to show, that
\medskip

{\sl a)}\ $U(y(k),\lambda , s)\, $ depends only on the   variable $k\, -$th

column of matrix $\lambda \in \Lambda_q\, ;$

{\sl b)}\ $U(y,\lambda , s)\,=\, \prod\limits_k\, U(y(k),\lambda , s)\, $
\medskip

We denote by $Y_q(\infty , q)\, $ a subset in $Y(\infty , q)\, ,$
that consists of matrices, which has the elements of first  $q\, $ rows
equal to zero. Let $L_q\, $
be a subgroup in $G_d(q,\infty)\in G \, ,$ that consists of matrices of
the form
$l_q(x)=\pmatrix{{({(l(x))}^{-1})}^{\prime}&0\cr
         0&{l(x)}\cr}\, ,$ where $\, l(x)= \pmatrix {{I_q}&0&0\cr
0&{I_q}&x\cr0& 0&{I}\cr}\, ,$ $x\, $
is an arbitrary local nonzero\ \  $q\times\infty \, -$ matrix.
 \medskip

Element $\lambda \in \Lambda_q \ $ we will write as $\lambda =
(\lambda_q,\lambda (q,\infty))\, ,$
where $\,\lambda_q \,  \ \mbox{is} \  \, q\times q\, -$ matrix, that
consists of first $q\, $
columns of matrix $\lambda\, .$
\medskip

If $y \in Y_q(\infty , q)\, ,$ then from the correlation
$\,l_q(x){\theta}^{(q)}(0,y)={\theta}^{(q)}(0,y)l_q(x) \, ,$ setting\newline
$\ {{\alpha}_{_A}}(\lambda,g) = \exp \lbrace
Tr \lbrack\, i\,  A\lambda(gg^{\ast}-1)\lambda^{\ast} \rbrack \rbrace\ $
and accounting (21) -- (22), we get
$$\ {{\alpha}_{_A}}((\lambda_q,\lambda (q,\infty))\, ,\, l(x))\,
U(y,(\lambda_q, \lambda (q,\infty)+\lambda_qx),s)=$$
\begin{equation}
 {{\alpha}_{_A}}((\lambda_q,\lambda (q,\infty)+2z^{-1}{(yh(s))}^t)\, ,\, l(x))
\, U(y,\lambda ,s)\, .
\end{equation}

Let

$ a(\lambda , x,y)=\exp\{\, -2i\, Tr[{\lambda}^*(q,\infty)Az^{-1
}{(yh(s))}^t+$

${({\lambda}^*(q,\infty)Az^{-1}{(yh(s))}^t)}^*+x^*{\lambda}_q^*Az^{-1}
(s){(yh(s))}^t+$

${(x^*{\lambda}_q^*Az^{-1}(s){(yh(s))}^t)}^*+2{(z^{-1}(s){(yh(s))}^t)}^*A
z^{-1}(s){(yh(s))}^t]\}\, .$

Since $\ {{\alpha}_{_A}}(\lambda,g) = \exp \lbrace
Tr \lbrack\, i\,  A\lambda(gg^{\ast}-1)\lambda^{\ast}
\rbrack \rbrace\ \, ,$
then using (24) we get
\begin{eqnarray}
&&a(\lambda,x,y)U(y,(\lambda_q, \lambda (q,\infty)+\lambda_qx),s)=\nonumber\\
&&a(\lambda,0,y)U(y,\lambda,s) .
\end{eqnarray}

From the properties ({\sl a}) -- ({\sl b}), accounting a form of the action
of operators $\, {\widehat\Pi}_s(g_q)\, ,\,
\newline{\widehat\Pi}_s({\theta}^{(q)}(0,y) )\, $
(see (20) and (22)) and the inclusion $y \, \in\,Y_q(\infty , q) \, ,$
it's easy to get that
$\, U(y,\lambda,s)\, $ does not depend on $\,\lambda_q\, .$
\medskip

Hence, accounting (25), we get that
$\, a(\lambda,0,y)U(y,\lambda,s) =
c(y,s)\, $ does not depend on $\, \lambda\, .$

Therefore,
\begin{eqnarray}
&&U(y,\lambda,s) =c(y,s) \exp\{\, -2i\,
Tr[{\lambda}^*(q,\infty)Az^{-1}{(yh(s))}^t+\nonumber\\
&&{({\lambda}^*(q,\infty)Az^{-1}{(yh(s))}^t)}^*+
2{(z^{-1}(s){(yh(s))}^t)}^*Az^{-1}(s){(yh(s))}^t]\}\nonumber
\end{eqnarray}
for all $\, y\in Y_q(\infty , q)\, $ and $\lambda \in \Lambda_q \, .$
\medskip

Taking into account this correlation it's  easy to check, that
${\cal U}(y,\lambda,s) =  U(y,\lambda,s) c^{-1}(y,s)\, $
is 1--cocycle of action of group $\,{\theta}^{(q)}(0,\,
Y_q(\infty , q))  \,$ to $\, \Lambda_q \, .$
Since $\,U \, $ is also 1--cocycle, and $\,c(y,s)\, $ does not depend on
$\,\lambda \, ,$
then $\,c(\cdot , s)\, $ is a multiplicative character of group $\,
Y_q(\infty , q)) \, .$
\medskip

If $\, g=\pmatrix{I_q&0\cr 0&g_1\cr}\, ,\, g_q=
\pmatrix{{(g^{-1})}^{\prime}&0&0&0\cr
0&I_q&0&0\cr 0&0&I_q&0\cr 0&0&0&g\cr}\, ,$ then from the correlation
\medskip

$\, {\widehat\Pi}_s(g_q) {\widehat\Pi}_s({\theta}^{(q)}(0,y) )
{({\widehat\Pi}_s(g_q))}^*={\widehat\Pi}_s({\theta}^{(q)}(0,{(g_1^{-1})}^ty)
 ) \, ,$
accounting (20),  (22) and using the simple calculations,  we get $\,
c(y,s)=c({(g_1^{-1})}^ty,s)  \, \,\, \forall\, g_1 \, .$
Therefore, $\, c\, \equiv 1\, .$ Correlation ({\sl i}) is proved.
\medskip

To prove ({\sl ii}) and ({\sl iii}) we get the correlation, that connects
$\, z(s)\, $ and $\, h(s)\,  (s\in S)\, $ (see (19)--(23)).
\medskip

Obviously the isometry $\, V_s \, (s\in S)\, $ from the proposition 6.7 maps
$\, \Pi_s(U(GK_q))  \, -$ fixed vector $\, {\eta}_q^{^{(U)}}(s)\, $ into the
vector $\, \xi_0 \in L^2(\Lambda_q,\nu_q)\, $  determined by the
function,
that  equals to unit  on $\, \Lambda_q\, .$
 \medskip

Let

$\, s_-^{(q)}(e_i)\, (\, s_+^{(q)}(e_i) \,)\ =\ \left\{
\begin{array}{rl}
s_-^{(q)}(e_i)\, (\, s_+^{(q)}(e_i) \,) ,&\mbox{if}\, |i|>q\\
e_i\ ,  &\mbox{if}\, |i|\leq q .
\end{array} \right. $
\medskip

Since $\, \Pi \, $ is an admissible representation, then it extends
by the continuity
to $\, s_-^{(q)} \ (\, s_+^{(q)} \,)  \ .$ Moreover, operators
$\, \Pi (s_-^{(q)})\ (\Pi (s_+^{(q)}))  \, $
we can approximate by the elements from $\, \Pi (U(GK_q)) \, .$ Hence
$\,\Pi_s(s_{\pm}^{(q)}){\eta}_q^{^{(U)}}(s)\, =\, {\eta}_q^{^{(U)}}(s)\,  \, $
(see lemma 6.3 and proposition 6.6). Therefore
$\, {\hat  \Pi}_s(s_{\pm}^{(q)}) \xi_0=\xi_0\, $
This condition and the correlation
$\, {(s_{\pm}^{(q)})}^{-1} {\theta}^{(q)}(x,0)  s_{\pm}^{(q)} \, =
\, {\theta}^{(q)}(0,\pm x)\  $
lead to the equality:
\medskip

$\, ({\hat  \Pi}_s( {\theta}^{(q)}(x,0) )\xi_0\, ,\,\xi_0\, )  \, =
\, ({\hat  \Pi}_s( {\theta}^{(q)}(0,\pm x) )\xi_0\, ,\,\xi_0\, ) \, .$

\medskip
We calculate the both sides of this equality using
({\sl (i)}) and the correlations
(21) -- (22), and we get
\begin{eqnarray}
&&\exp\bigg{\{}-{1\over 4}Tr(x^*xz(s)z^*(s))\bigg{\}}=\nonumber \\
&&\exp\{-\, Tr[x^*x(4h(s)({{(z^{-1})}^*A^2z^{-1})}^th^*(s)+h(s)
({{(z^{-1})}^*z^{-1})}^th^*(s))]\}\nonumber
\end{eqnarray}

Therefore,
\begin{eqnarray}
z(s)z^*(s)=4h(s){(z^{-1}(s))}^t{(1+4A^2)}^t{({(z^{-1}(s))}^t)}^*h^*(s)\ .&&
\end{eqnarray}
Since $\, {\it rang }\ z(s)=q \,$ for $\mu\, -$a.a. $s\in S \,$
(see proposition 6.1), then from
(26) we get the property ({\sl ii}) .
\medskip

Let's pass to the proof of ({\sl iii}) .

At first we notice that from ({\sl ii}) it follows the correctness of
the definition of
matrix $\, R(s)={1\over 4}z^{t}(s)h^{-1}(s)z(s) \, .$
\medskip

Put
\begin{eqnarray}
T(b,\lambda ,s)=\exp\{ -i\, \Re Tr[R(s)\lambda b{\lambda}^t]\} \ .\nonumber
\end{eqnarray}
Correlation $\, {\theta}^{(q)}(0,y) {\gamma }_{_{\scriptstyle u}}^{(q)}(-b)
{\theta}^{(q)}(0,-y) ={\gamma }_{_{\scriptstyle u}}^{(q)}(-b)
{\delta}^{(q)}({(by)}^ty) {\theta}^{(q)}(by,0) \, $
 (see (8)) leads to the equality:
\begin{eqnarray}
M(-b,\lambda +2z^{-1}{(yh(s))}^t,s)=M(-b,\lambda ,s)\exp\{i\,
\Re Tr\{ [{(by)}^t yh(s)+z(s)\lambda by]\} \ .\nonumber
\end{eqnarray}
Hence, accounting the definition of $\, T(b,\lambda ,s) \, ,$ we get
\begin{eqnarray}
\, M(-b,\lambda +2z^{-1}{(yh(s))}^t,s)T(b,\lambda +2z^{-1}{(yh(s))}^t ,s)=
M(-b,\lambda ,s)T(b,\lambda ,s)\ .\nonumber
\end{eqnarray}
Since $y\, $ is an arbitrary, Ánd $\, {\it rang}\ h(s)= {\it rang}\ z(s)=q \, ,$
then
$\, M(-b,\lambda ,s)T(b,\lambda ,s)=c(b,s) \, $does not depend on $\,
\lambda \, .$
 Hence, accounting the definition of $\, T(b,\lambda ,s)\, $ and correlation
$\, g_q {\gamma }_{_{\scriptstyle u}}^{(q)}(b) {(g_q)}^{-1}=
{\gamma }_{_{\scriptstyle u}}^{(q)}(gbg^t) \,  $
(see (8)), we get $\, c(b,s)=c(gbg^t,s)\ \ \forall g   \, .$ Therefore,
$c(b,s)=1\, $ for $\,\mu \, -$ a.a. $\, s\in S\, .$ Lemma 6.9 is proved.
\medskip

{\bf Lemma 6.10}\ \ {\it Let $\, R(s)={1\over 4}z^t(s)h^{-1}(s)z(s)\, ,\,
\newline {|R(s)R^*(s)|}^{1\over 2}u(s)=R(s)\, $ is a polar decomposition of
$\, R(s) \, .$
\newline ôhen for $\,\mu-$a.a. $\ s\in S\, $ the following correlations are
true:

{\sl i)}\ $\, I+4{(A^t)}^2=4R(s)R^*(s)\ ;$

{\sl ii)}\ $\, -A=u^*(s)A^tu(s)\ .$}
\medskip

{\bf Proof} \ Let $\ g_q  =
\pmatrix{{(g^{-1})}^t &0&0&0 \cr 0& I_q &0&0\cr 0&0& I_q &0\cr 0&0&0&g\cr}
\ ,\
\xi_0 \ $ is a vector from $\, L^2(\Lambda_q,\nu_q) \, $ defined by the
function,
that equals to unit. Suppose, that matrix $\, g \, $ has a form
$\, \pmatrix{{g(n)}&0\cr 0&I\cr}\ ,$ where $\, g(n) \, \  \mbox{is} \ \ \,
n\times n\, -$ matrix,
for which $\, (g(n)g^*(n)-I_n ) \, $ is invertible. Identify  $\,
C^{qn} \, $
with a set of all complex $\, n \times q \,-$matrices .
\medskip

By (20)

$\, ({\widehat\Pi}_s(g_q)\xi_0 )(\lambda)=
\exp {\bigg [}-{1\over 2}Tr((1-2\,i\, A)\lambda (gg^*-1){\lambda}^*)
{\bigg ]}\, .$
\medskip

Hence and from (21) using the ordinary calculations we get
\begin{eqnarray}
&&({\widehat\Pi}_s(g_q)\xi_0 )(\lambda)=c(g,s)\int\limits_{C^{qn}}
\exp[-{1\over 2}Tr(z(s){1-2\,i\, A)}^{-1}\cdot\nonumber \\
&&\cdot z^*(s)x_n^*{(g(n)g^*(n)-I_n )}^{-1}x_n)]({\widehat\Pi}_s
({\theta}^{(q)}(x,0) )\xi_0 )(\lambda)\, dx_n\, ,\ \ \ \ \ \ \ \ \ \
\end{eqnarray}
where $x_n\in C^{qn}\, ,\, x=\pmatrix{x_n\cr 0\cr} \, ,\, c(g,s)\, $
depends only on $\, g\,$ and $\, s \, .$
\medskip

From the definition of $\, s_{\pm}^{(q)} \,$ and from (27) the next chain
of
equalities follows:

\begin{eqnarray}
&&({({\widehat\Pi}_s(s_{\pm}^{(q)}))}^{-1}{\widehat\Pi}_s(g_q)
{\widehat\Pi}_s(s_{\pm}^{(q)})\xi_0)(\lambda)=
({\widehat\Pi}_s({(g_q^{-1})}^t)\xi_0)(\lambda)= \nonumber \\
&&\exp {\bigg {\{}} -{1\over 2}Tr{\bigg [}(1-2\, i\, A)
\lambda ({(g^{-1})}^t {\bar g}^{-1}-1){\lambda }^*{\bigg ]}{\bigg {\}}}=
c(g,s)\cdot\nonumber \\
&&\cdot\int\limits_{C^{qn}}\exp {\bigg [}-{1\over 2}
Tr(z(s){1-2\,i\, A)}^{-1} z^*(s)x_n^*{(g(n)g^*(n)-I_n )}^{-1}x_n)
\bigg {]}\cdot \nonumber \\
&&\cdot ({\widehat\Pi}_s({\theta}^{(q)}(0,\pm x) )\xi_0 )
(\lambda)\, dx_n\, .\nonumber
\end{eqnarray}
Next, accounting (22) and  a statement of lemma 6.9({\sl i}), we get:
\begin{eqnarray}
&&\exp {\bigg {\{}} -{1\over 2}Tr{\bigg [}(1-2\, i\, A)
\lambda ({(g^{-1})}^t {\bar g}^{-1}-1){\lambda }^*
{\bigg ]}{\bigg {\}}}= c(g,s)\cdot\nonumber \\
&&\cdot\int\limits_{C^{qn}}\exp {\bigg {\{} }-{1\over 2}
Tr{\bigg [}z(s){(1-2\,i\, A)}^{-1} z^*(s)x_n^*{(g(n)g^*(n)-I_n )}^{-1}
x_n-           \nonumber \\
&&-2(1-2\, i\, A)z^{-1}(s)h(s)x^t{\lambda}^*-2h^*(s){(z^*(s))}^{-1}
(1-2\, i\, A)\lambda {\bar x}+                    \nonumber \\
&&+4{\bar x}h^*(s){(z^*(s))}^{-1}(1-2\, i\, A)z^{-1}(s)h(s)x^t
{\bigg ]}{\bigg {\}} }\, dx           \, .                                    \nonumber
\end{eqnarray}
The calculation of integral in the right side leads to the correlation
\begin{eqnarray}
\exp {\bigg {\{ }} - {1\over 2}Tr{\bigg [ }(1-2\, i\, A)
\lambda ({(g^{-1})}^t {\bar g}^{-1}-1){\lambda }^*{\bigg ] }
{\bigg {\}}}  =\phantom{ 00000000000000000000 } && \nonumber \\
\exp {\bigg{\{ } }2Tr{\bigg [ }(1-2\, i\, A){{\bigg (}
z^*(s){(h^*(s))}^{-1}{\bar z }(s){(1-2\, i\,  A^t  )}^{-1}
z^t(s)h^{-1}(s)z(s)\otimes }\phantom{0 } && \\
{\otimes {{\bigg (} {(g(n)g^*(n))}^t-1{\bigg ) }}^{-1}+4(1-2\,
i\, A){\bigg )}}^{-1}\{ (1-2\, i\, A )\lambda (n)\}
{\lambda }^*(n){\bigg ] }{\bigg {\} }} .\phantom{z^t(s)} &&  \ \nonumber
\end{eqnarray}
Here $\, \lambda (n) \,$ consists of first $\, n\,$ columns of matrix $\,
\lambda  \, ,$
the action of the operator $\, a\otimes b \, $ on $\, \lambda \in
\Lambda_q \, $ is defined
by: $\, (a\otimes b)\{\lambda \}=a\lambda b \, .$
\medskip

If we use the ordinary calculation, we'll see that  (28) is true if and only
if
 $\, 4(1-2\, i\, A) =\newline z^*(s){(h^*(s))}^{-1}{\bar z }(s)
{(1-2\, i\,  A^t  )}^{-1}z^t(s)h^{-1}(s)z(s) \, .$
Therefore,
\begin{eqnarray}
&& 1+4{(A^t)}^2=4R(s)R^*(s) \nonumber \\
&& -4A=R^*(s){(1+4{(A^t)}^2)}^{-1}A^tR(s)\ .
\end{eqnarray}
If $\, {|R(s)R^*(s)|}^{1\over 2}u(s)=R(s)\, $ is a polar decomposition of
$\, R(s) \, ,$
then from (29) it follows that  $\, -A=u^*(s)A^tu(s)\ .$ Lemma 6.10 is
proved.

\bigskip
For $\, \mu \, -$a.a. $\, s\in S \, $  representation $\,
{\widehat\Pi}_s \, $
of group $\, GK_q  \, $ is defined by the set of parameters
$\, \{A\, ,\, z(s)\, ,\, h(s)\}\, ,$   \,
that was introduced in propositions 6.7 -- 6.8. If $\, \varepsilon\,
$is a non negative\newline number,
 then we denote by $\, \vartheta  ( \varepsilon )\, ,\, \varsigma
( \varepsilon ) \, $
the   $\, 2\times 2\, -$matrices of the form $\, \varsigma
( \varepsilon) =\pmatrix{0&{i\varepsilon }\cr
{-i\varepsilon}&0\cr} \,  , \,
\vartheta ( \varepsilon)=\pmatrix{{\varepsilon}&0\cr 0&{-\varepsilon}\cr}
\, .   $

From the propositions 6.9 -- 6.10 using the standard methods of elementary
theory of matrices
we can get the following statement.
\medskip

{\bf Proposition 6.11.}\ \ {\it There exists $\, \mu \, -$measurable mapping
$\, w(\cdot ) \,$ from $\, S \,$ into  $\, U(q) \,$ with the following
properties:

{\sl i)}\ if $\, G=Sp(2\infty ) \, ,$ then

$w(s)Aw^*(s)=d_{sp}(A)=\pmatrix{{ \varsigma (\lambda_1)}& 0_2 & \ldots & 0_2
 & 0\cr
                               0_2 &{ \varsigma (\lambda_2)}& \ldots & 0_2&
0\cr
                               \vdots & \vdots & \ddots & \vdots & \vdots &
\cr
                                0_2 & 0_2 & \ldots & { \varsigma (\lambda_k)} &
 0\cr
                                0&0&0&0&0_p\cr }  \, ,$

where $\, 0_p  \, $ is zero $\, p\times p\, -$matrix, $\, \lambda_j >0   \,
$ for

$ j=1,2,\ldots ,k\, $ ($\, k \,$ may be equals to zero  )

and $\, 2k+p=q \, ;$

$\, {\bar w}(s)u(s)w^*(s)=I_q \, $ (see lemma 6.10) ;
\medskip

{\sl ii)}\ \ if $\, G=O(2\infty ) \, ,$ then

$w(s)Aw^*(s)=d_o(A)=\pmatrix{{ \vartheta (\lambda_1)}& 0_2 & \ldots &
0_2 & 0\cr
                               0_2 &{ \vartheta (\lambda_2)}& \ldots &
0_2& 0\cr
                               \vdots & \vdots & \ddots & \vdots &
\vdots & \cr
                                0_2 & 0_2 & \ldots & { \vartheta (\lambda_k)} & 0\cr
                                0&0&0&0&0_{2p}\cr }  \, ,$

where $\, 2p+2k=q \, ,$ Ánd the rest properties of the parameters

of matrix $\, d_o(A) \,$ are the same as in ({\sl i});

$\, u_o={\bar w}(s)u(s)w^*(s)=\pmatrix{{ \varsigma (-i)}& 0_2 & \ldots &
 0_2 \cr
                               0_2 &{ \varsigma (-i)}& \ldots & 0_2\cr
                               \vdots & \vdots & \ddots & \vdots &  \cr
                    0_2 & 0_2 & \ldots & { \varsigma (-i)} \cr }  \, ;$

\medskip

{\sl iii)}\ \ $\, R_c= {\bar w}(s)R(s)w^*(s)={1\over 2}
\sqrt{\mathstrut {1+4w(s)A^2w^*(s)}}\, {\bar w}(s)u(s)w^*(s)\, $

does not depend on $\, s\, .$}
\bigskip

Consider the unitary operator $\, W\, ,$ that acts in
$\, L^2(S,\mu)\otimes L^2(\Lambda_q ,\nu_q)\, $ by
$$(W\eta )(s,\lambda )=\eta (s,w^*(s)\lambda)\, ,\ \
{\mbox{where}}\ \ \eta\in L^2(S,\mu)\otimes L^2(\Lambda_q ,\nu_q)\, .$$
\bigskip

The next statement follows from the
propositions 6.6 -- 6.8, 6.11 and has the same
meaning for description of the representations of
group $\, G \, ,$ as the lemma 5.5
in a case of $\, GL(\infty ) \, .$
\bigskip

{\bf Proposition 6.12.}\ \ {\it We denote by $\, V  \, $ the isometry of
a Hilbert space $\newline H_q^{^{(U)}}=\int\limits_S^{\oplus}H_q^{^{(U)}}(s)\, d
\mu (s)\, $
 to $\, L^2(S,\mu)\otimes L^2(\Lambda_q ,\nu_q)\, ,$ that defined by
 $\, (V\eta)(s)=V_s\eta (s) \, ,$ where $\, \eta (s)\in
H_q^{^{(U)}}(s) \, ,\,
V_s\eta (s)\in L^2(\Lambda_q,\nu_q)\, $ for $\, \mu\, -$a.a. $\, s\in S \, $
(see propositions 6.6 -- 6.7). For $\, g \in GK_q  \,$ we put
\begin{eqnarray}
&&\ddot{\Pi}(g)=WV\Pi (g)V^*W^*\ ,\ \ \zeta (s)=z(s)w^*(s)\ . \nonumber
\end{eqnarray}
\def\ort {\scriptstyle d_o(A)}  \def\si {\scriptstyle d_{sp}(A)}

ôhen the action of  operators $\ \ddot{\Pi}(g)\ $ in $\, L^2(S,\mu)
\otimes L^2(\Lambda_q ,\nu_q) \,$
is defined by  correlations:

\begin{equation} ({\ddot\Pi}({\theta}^{(q)}(x,0) )\eta )(s,\lambda)\, =\,
\exp \lbrack i  \Re  Tr (\zeta (s) \lambda h )
\rbrack \eta (s,\lambda )\, ;
\end{equation}
\begin{center}$({\ddot\Pi }(g_q)\eta )(s,\lambda )=$
\end{center}\begin{equation}\left\{
\begin{array}{rl}
{\hat{\alpha }_{_{{\scriptstyle d_{sp}(A)} }}}(\lambda ,g)
\eta (s,\lambda g), \mbox{if   } g_q \in Sp(2\infty)&\\
{\hat{\alpha }_{_{\scriptstyle d_o(A) }}}(\lambda ,g)\eta (s,\lambda g),      \mbox{if   }g_q \in O(2\infty) \ ;  &
\end{array} \right.
\end{equation}
\begin{center} $( {\ddot\Pi}({\gamma }_{_{\scriptstyle u}}^{(q)}(b) )\eta )
(s,\lambda)=$\end{center}
\begin{equation}\left\{
\begin{array}{rl}
\exp{\bigg{\{}}-{{i }\over 2 } Tr[\sqrt{\mathstrut {1+4{(d_{sp}(A))}^2}}
\, \lambda b {\lambda}^t ]{\bigg {\}}}   ,\mbox{if   } b=b^t \\
\exp{\bigg{\{}}-{{i }\over 2 } Tr[\sqrt
{\mathstrut {1+4{(d_{o}(A))}^2}}\, u_o\lambda b {\lambda}^t ]{\bigg {\}}},
\mbox{if } b=-b^t \
\end{array} \right.
\end{equation}
(see proposition 6.11({\sl ii}) and lemma 6.9({\sl iii}));
$$({\ddot\Pi}({\theta}^{(q)}(0,y) )\eta )(s,\lambda) =
U(y,\lambda , s)\,\cdot $$\begin{equation}
\cdot {\bigg{[}{{d\nu_q(\lambda + 2{\zeta}^{-1}(s)
(y{h}(s))^t)}\over{d\nu_q(\lambda) }}\bigg{]}}^{1\over{2} }
\eta )(s\, ,\,\lambda + 2{\zeta}^{-1}(s)(y{h}(s))^t)\ .
\end{equation}

The form of the unitary cocycle from the last correlation is given in the
condition of lemma 6.9 (see (22)).
Moreover in the corresponding  expression we should change
$\, z(s)\,$ into $\, \zeta (s) \, ,$
and $\, A\,$-- into $\, w(s)Aw^*(s) \, $ (see proposition 6.11).}

\bigskip

The advantage of the realization of $\, {\ddot\Pi} \,$ by correlations
(20) -- (23)
is that the form of the action of operators $\, {\ddot\Pi}(g_q)  \,$
and $\, {\ddot\Pi}({\gamma }_{_{\scriptstyle u}}^{(q)}(b) \, $
does not depend on $\, s\, .$
\medskip

Since $\, {\ddot\Pi} \, $ is an admissible representation, then
it extends by continuity
to the group, generated by $\, GK_q \,$  and $\, s_-^{(q)} \,$ for
$\, G=Sp(2\infty)  \,$ or $\, GK_q \,$  and $\, s_+^{(q)} \,$ for
$\, G=O(2\infty)  \, .$

Change if we need  $\, \Pi (U(GK_q)) \, -$
fixed vector $\, {\eta }_q^{^{U}} \, $ (see lemma 6.3 and proposition 6.6 )
to $\, c{\eta }_q^{^{U}}  \, ,$ where $\, c\, $  is an operator from the
center of  $w^*-$ algebra $\, {(\Pi (GK_q))}^{\prime\prime}\, ,$
we may take    $\, WV{\eta }_q^{^{U}} =\xi_0 \,$ (see proposition 6.12).
Vector $\, \xi_0\,$ is determined by the function on $\, S\times
\Lambda_q \, ,$
that equals to the unit. Therefore, $\, {\ddot\Pi}(s_{\pm}^{(q)})\xi_0=
\xi_0 \, .$
Hence and from the correlations (20), (22) we get
\newline
\medskip
$\, {({\ddot\Pi}(s_{\pm}^{(q)}))}^{-1}{\ddot\Pi}({\theta}^{(q)}(x,0)\xi_0=
 {\ddot\Pi}({\theta}^{(q)}(0,\pm x)\xi_0\, .$

Thus
\begin{equation}
\begin{array}{rl}
&\exp \{ i  \Re  Tr (\zeta (s) \lambda h ) \}
\eta (s,\lambda ) \stackrel{ {{\ddot\Pi}(s_{\pm}^{(q)})}^{-1}}
{\longrightarrow} \\
&{\longrightarrow}\exp\{2i Tr[\pm {\lambda }^*
d_{\sharp}(A){\zeta}^{-1}(s)h^t(s){({\zeta}^t(s))}^{-1}y^t\pm \\
&\pm {\bar y}{{\bar{\zeta (s)}}^{-1}}{\bar h (s) }{\zeta }^{-1}(s)
d_{\sharp }(A)\lambda +\\
&+2 {\bar y}{{\bar{\zeta (s)}}^{-1}}{\bar h (s) }{\zeta }^{-1}(s)
d_{\sharp }(A){\zeta}^{-1}(s)h^t(s){({\zeta}^{-1})}^ty^t]\}\cdot \\
&\cdot {\bigg [}{{\displaystyle d\nu_q(\lambda +2{\zeta}^{-1}(s)h^t(s)
{({\zeta }^{-1}(s))}^ty^t)}\over {\displaystyle d\nu_q(\lambda )}}
{\bigg ]}^{1\over 2}\, ,
\end{array}
\end{equation}
where $\,  \sharp \,$ means one of the indexes $\ o\ $ or $\, sp \, .$
\vskip 6pt

From this correlations we get an important

{\bf Remark 6.13.} \ {\it  Since $\, \zeta (s)=z(s)w^*(s) \, ,$
then by the statements 6.10 -- 6.12 the next correlation is true

$$R=\bar w (s)R(s)w^*(s)={1\over 4}{\zeta}^t(s)h^{-1}(s){\zeta}(s)\, ,$$
and that's why the form of the action of the operator $\,
{\ddot\Pi}(s_{\pm}^{(q)}) \, $
 does not depend on $\, s\in S \, .$}
\bigskip

We denote by $\, \widetilde G(q,\infty) \,$ a group generated by all
the elements
$\, g_q \, ,\, {\gamma }_{_{\scriptstyle u}}^{(q)}(b)\, ,\, s_{\pm}^{(q)}\, .    $
Obviously $\,  G(q,\infty) \subset  \widetilde G(q,\infty) \, .$ Therefore,
for $\, g\in G(q,\infty) \, \ {\ddot\Pi}(g)=I_{_{L^2(S,\mu)}}\otimes
\widetilde\Pi (g) \, ,$
where $\,\widetilde\Pi (g) \, \,$ is a unitary operator in $\,
L^2(\Lambda_q,\nu_q)  \, $
(see statements 6.12 -- 6.13). We denote by $\, {\ell}_q \, $ a natural
isomorphism of groups $\, G(q,\infty)  \, $ and $\, G \, $\, $ ({\ell}_q:\,
g\in  G(q,\infty)\,\longrightarrow {\sigma}_q^*g{\sigma}_q \in G\,)\, .$
 \bigskip

The total of our reasoning is the following statement.
\medskip

{\bf Theorem 6.14. } {\it \  Let $\, \Pi  \, $ be the same as in propositions
6.1 -- 6.4, $\, \xi \,$ is a cyclic vector for $\, \Pi \, ,$ $ \, \xi_q=
\Pi (\sigma_q ) \, $
and $\, \Pi_A \,$ is a representation of group $\, G \, $ built in \S 1
(see propositions 1.2~-- 1.3).
ôhen the restriction of $\, \Pi \, $ to $\, G(q,\infty)\, ,$ that acts in
$\, [\Pi (G(q,\infty))\xi_q]\, ,$
is multiple by the irreducible component of the representation, that defined
by the chain of mappings:
$$g\in  G(q,\infty)\,\longrightarrow {\ell}_q(g)\in G\longrightarrow
\Pi_A ( {\ell}_q(g)) \ .$$}
\medskip

{\bf To prove  } this theorem it suffices to notice that $\,
\Pi_A\circ\ell_q  \,$
is naturally unitary equivalent to $\,\widetilde\Pi  \,  .$

 \bigskip

To formulate the main classification result
let us notice the important moments of our reasoning.
\medskip

We have supposed , that the admissible factor--representation $\, \Pi\,$
of group $\, G\,$
is cyclic with the respect to the unit vector $\, \xi\in H_{\Pi} \,$
(see remark 6.2).
Then we've considered a restriction of $\,\Pi \,$ to the subgroup
$\, G(q,\infty)\, ,$
that acts in $\, [\Pi ( G(q,\infty))\Pi  (\sigma_q)\xi] \, .$
There was shown that ,
it is a subrepresentation of the spherical representations of larger group
$\, GK_q\, $ (see proposition 6.4), for which in the proposition 6.12
explicit realization was given.
Moreover the  form of the operators $\, \ddot\Pi (g) \,$ (see (31), (32)
and (34)) for
$\, g \in  G(q,\infty) \,$ does not depend on $\, s\in S \, $
(see also remark 6.13).\
At last in theorem 6.14 we've got the form of  representation
$\, (\Pi\, ,\, G(q,\infty)\, ,\,[\Pi ( G(q,\infty))\Pi  (\sigma_q)\xi])\, .$
Theorem 0.5 found by G.I. Ol'shansky  let us turn
to the representation $\, \Pi \,$ of group $\, G\,$ .
\medskip

Because of this fact $\, \Pi\,$ extends by continuity to $\, \sigma_q\, $ .
Therefore the next correlation is true:
\begin{eqnarray}
{(\Pi (g)\xi ,\xi)}_{_{H_{\Pi}}}={(\Pi (\sigma_q  g{\sigma_q}^*+I_q(0) )
\Pi (\sigma_q)\xi , \Pi (\sigma_q)\xi )}_{_{H_{\Pi}}}\, ,
\end{eqnarray}
where $  I_q(0)=\pmatrix{0&0&0&0\cr
                       0&I_q&0&0\cr
                       0&0&I_q&0\cr
                       0&0&0&0\cr}\ ,$
Ánd $\, \ \ \sigma_q  g{\sigma_q}^*+I_q(0) \,$ belongs to
$\, G(q,\infty)\, .$
\medskip

Next, by theorem 6.14 in  $L_q^{^A}$  (see \S 1) ,
there exists  vector $\, f_{\xi }\, ,$ for which
\newline
\medskip
${(\Pi (\sigma_q  g{\sigma_q}^*+I_q(0) )\Pi (\sigma_q)\xi ,
\Pi (\sigma_q)\xi )}_{_{H_{\Pi}}}={(\Pi_A\circ\ell_q
(\sigma_q  g{\sigma_q}^*+I_q(0)) f_{\xi }\, ,\,  f_{\xi })}_{_{L_q^{^A}}}\ .$
\newline
Hence and from  (35) we get
$ \ {(\Pi (g)\xi ,\xi)}_{_{H_{\Pi}}}={(\Pi_A (g)f_{\xi },
f_{\xi })}_{_{L_q^{^A}}}\ .$
\medskip

This correlation and theorem 1.4 lead us to the main result of this chapter.
\medskip

{\bf Theorem 6.15.}{\it  \ Let $\, \Pi\, $ be an admissible
factor--representation
of group $G$\ \ $(G=Sp(2\infty)\mbox{or}\newline \ O(2\infty))\ .$
ôhen there exist:
a natural number $\, q\, ,\ q\times q\, -   $
self-adjoint matrix $\, A\, ,$
that satisfies to the corresponding conditions of
propositions 1.2--1.3, \ an admissible representation
$\, \rho \,$ of group $\ O(A,q)\subset U(q) \ $
for $\  G=Sp(2\infty) \ $ or
$\, Sp(A,m)\subset U(q)\, $  for $\  G=O(2\infty) \ $ (see \S 1)
such that
$\  \Pi\  $ is multiple by restriction of $\, \Pi_A\, $ to $\,
P_k^{\rho}  L_q^{^A}  \, $\ (see theorem 1.4).}

\end{document}